\documentclass{compositio}

\usepackage{kbordermatrix}		
\usepackage[T1]{fontenc}
\usepackage{ae,aecompl}     		
\usepackage{graphicx}
\usepackage[english]{babel}
\usepackage{marvosym}
\usepackage[tight,nice]{units}
\usepackage[bottom,multiple]{footmisc}  
\usepackage[mathletters]{ucs}		
\usepackage[utf8x]{inputenc}		
\usepackage{mabk} 			
\usepackage{url}			
\usepackage{amsmath}
\usepackage{amssymb}
\usepackage{amsfonts}
\usepackage{amstext}
\usepackage{mathrsfs}
\usepackage{bbm}
\usepackage[swap,section]{mtheorems}
\usepackage{dcpic,pictex} 		
\usepackage{hyperref} 			
\usepackage[numbers,square]{natbib	}

\newcommand{\hldef}[1]{\textbf{#1}}

\newcommand{\quadtextf}{\quad \text{ for }}
\newcommand{\tr}{\text{tr}}

\newcommand{\Inv}{\text{Inv}}

\newcommand{\im}{\mathbf{i}}

\newcommand{\pa}{\partial}
\let\f=\frac
\let\set=\mathbb
\newcommand{\norm}[1]{\Arrowvert #1 \Arrowvert}
\newcommand{\abs}[1]{| #1 |}

\newcommand{\aequiv}{\Leftrightarrow}

\newcommand{\dom}[1]{\text{dom}(#1)\xspace}

\newcommand{\ran}[1]{\text{ran}(#1)\xspace}

\newcommand{\Sec}[1]{\Gamma(#1)}
\newcommand{\Secc}[1]{\Gamma_c(#1)}

\newcommand{\tensor}{\otimes}

\newcommand{\isom}{\cong}

\newcommand{\End}{\text{End}}

\RequirePackage{relsize}

\newcommand{\mat}[1]{\text{M}(#1,\set R)}
\newcommand{\Mat}[1]{\text{M}(#1,\set C)}

\newcommand{\Rpq}{{\set R^{p,q}}}

\newcommand{\pr}{\text{pr}}

\newcommand{\Cl}{\text{Cl}}
\newcommand{\Clc}{\Cl^{\set C}}

\newcommand{\Pin}{\text{Pin}}
\newcommand{\GL}{\text{GL}}

\newcommand{\Spin}{\text{Spin}}
\newcommand{\SO}{\text{SO}}

\newcommand{\supp}{\text{supp}\,}
\newcommand{\dd}[1]{\frac{\pa}{\pa #1}}

\newcommand{\pq}{{p,q}}

\newcommand{\dvol}{dV}
\newcommand{\dvolMg}{dM_g}

\newcommand{\myispr}[1]{[#1]}		
\newcommand{\myspr}[1]{(#1)}		

\newcommand{\sprpq}[1]{\langle #1\rangle_\pq}
\newcommand{\sprRn}[1]{\myspr{#1}_n} 
\newcommand{\sprRndummy}{\sprRn{\cdot,\cdot}} 
\newcommand{\sprpqdummy}{\sprpq{\cdot,\cdot}} 

\newcommand{\sprD}[2][\pq]{\myspr{#2}_{\Delta_{#1}}} 
\newcommand{\sprDdummy}[1][\pq]{\sprD[#1]{\cdot,\cdot}} 
\newcommand{\isprD}[2][\pq]{\myispr{#2}_{\Delta_{#1}}}
\newcommand{\isprDdummy}{\isprD{\cdot,\cdot}}

\newcommand{\sprB}[2][\xi]{\langle#2\rangle^S_{#1}}	
\newcommand{\sprBdummy}[1][\xi]{\sprB[#1]{\cdot,\cdot}}
\newcommand{\isprB}[1]{\myispr{#1}^S}
\newcommand{\isprBdummy}{\isprB{\cdot,\cdot}}
\newcommand{\sprBF}[1]{\langle#1\rangle^{S_F}}
\newcommand{\sprBFdummy}{\sprBF{\cdot,\cdot}}

\newcommand{\sprL}[2][\xi]{\myspr{#2}_{#1}}		
\newcommand{\isprL}[1]{\myispr{#1}}		
\newcommand{\sprLdummy}{\sprL{\cdot,\cdot}}
\newcommand{\isprLdummy}{\isprL{\cdot,\cdot}}

\let\spr\myspr 
\newcommand{\sprdummy}{\myspr{\cdot,\cdot}} 
\newcommand{\isprdummy}{\myispr{\cdot,\cdot}}

\newcommand{\normEucl}[1]{\norm{#1}_n}
\newcommand{\normx}[1]{\norm{#1}_\xi}	
\newcommand{\normdummy}{\norm{\cdot}}	
\newcommand{\normEucldummy}{\normEucl{\cdot}}

\newcommand{\specap}{\sigma_{ap}}
\newcommand{\specc}{\sigma_{c}}
\newcommand{\specr}{\sigma_{r}}
\newcommand{\specp}{\sigma_{p}}

\let\ker\kerold
\newcommand{\ker}[1]{\text{ker}(#1)}

\let\mcite\citep
\newcommand{\mcitep}[2]{\citep[p.\,#1]{#2}} 		
\newcommand{\mcitesatz}[2]{\citep[Satz #1]{#2}} 	
\newcommand{\mcitethm}[2]{\citep[thm.\,#1]{#2}} 	
\newcommand{\mciteprop}[2]{\citep[prop.\,#1]{#2}} 	
\newcommand{\mcitelem}[2]{\citep[lem.\,#1]{#2}} 	

\newcommand{\abbildungsdef}[5]{\begin{array}{llll} #1:  &#2&\to &#3\\%
						  &#4&\mapsto &#5\\ \end{array}}
\newcommand{\FT}{\mathcal{F}}
\newcommand{\id}{\text{id}}
\newcommand{\Id}{\text{I}}

\newcommand{\esssup}{\text{ess sup}\,}

\newcommand{\linearspan}[1]{\text{span}(#1)}

\renewcommand{\ale}{almost everywhere\xspace}
\providecommand{\proofname}{Proof}

\makeatletter
\newcommand\xleftrightarrow[2][]{%
  \ext@arrow 9999{\longleftrightarrowfill@}{#1}{#2}}
\newcommand\longleftrightarrowfill@{%
  \arrowfill@\leftarrow\relbar\rightarrow}
\makeatother

\usepackage{enumitem}
\setlist[enumerate]{label*=(\alph*),ref=(\alph*)}

\begin{document}
\title{Remarks on the spectrum of the Dirac operator of pseudo-Riemannian spin manifolds}
\author{Momsen Reincke}
\email{mail@momsenreincke.de}
\date{}
\begin{abstract}
	We study the spectrum of the Dirac operator $D$ on pseudo-Riemannian spin manifolds of signature $(\pq)$, considered as an unbounded operator in the Hilbert space $L^2_\xi(S)$. The definition of $L^2_\xi(S)$ involves the choice of a $p$-dimensional time-like subbundle $\xi\subset TM$.
We establish a sufficient criterion for the spectra of $D$ induced by two maximal time-like subbundles $\xi_1,\xi_2\subset TM$ to be equal.
If the base manifold $M$ is compact, the spectrum does not depend on $\xi$ at all.
We then proceed by explicitely computing the full spectrum of $D$ for $\set R^\pq$, the flat torus $\set T^\pq$ and products of the form $\set T^{1,1}\times F$ with $F$ being an arbitrary compact, even-dimensional Riemannian spin manifold.
\end{abstract}

\maketitle

\tableofcontents
\section{Introduction}

We investigate the spectrum of the Dirac operator $D$ on pseudo-Riemannian\footnote{The adjective \emph{pseudo-Riemannian} always refers to an indefinite metric in the tangent bundle, whereas the adjective \emph{semi-Riemannian} includes indefinite as well as definite metrics.} spin mani\-folds. 
The Dirac operator $D$ is a first-order differential operator acting on the smooth sections $\Sec{S}$ of the spinor bundle $S$.
In contrast to the Riemannian case, there is no canonical pre-Hilbert space structure (\ie a canonical positive-definite inner product) defined on the space $\Secc{S}$ of compactly supported sections of $S$. Therefore, there is no canonical Hilbert space in which $D$ acts as an unbounded Hilbert space operator.
Following \mcite{Baum1981}, one can define a positive-definite bundle metric in $S$ by fixing a maximal time-like subbundle $\xi$ of the tangent bundle $TM$. By means of this construction one defines a positive-definite inner product on $\Secc{S}$ and a Hilbert space $L^2_\xi(S)$. 
The Dirac operator then is an unbounded operator in this Hilbert space.
$L^2_\xi(S)$ additionally carries a Krein space structure (\cf \mcite{Bognar1974,Baum1981}).
Using Krein space arguments it can be shown that the spectrum of the Dirac operator has some symmetries \wrt the real and the imaginary axis (\mcite{Baum1981, Baum1994}).

How does the spectrum of $D$ depend on $\xi$?
This is an open question we will only partially answer.
After fixing some notation in section \ref{sec:notation}, we introduce the concept of quasi-isometric bundle metrics in a general vector bundle $E$ in section \ref{sec:QuasiIsometricBundleMetrics}.
Quasi-isometric bundle metrics introduce canonically isomorphic $L^2$-spaces $L^2(E)$.
Using the notion of quasi-isometric bundle metrics, we then study the dependence of the spectrum of the Dirac operator on the choice of $\xi$. 
More specifically, in theorem \ref{thm:DiracSpectrumsCoincide}, we will show that if the associated Riemannian metrics $r_{g,\xi_1}, r_{g,\xi_2}$ induced by two maximal time-like subbundles $\xi_1,\xi_2\subset TM$ are quasi-isometric, then the spaces $L^2_{\xi_1}(S)$ and $L^2_{\xi_2}(S)$ are isomorphic in a way compatible with $D$ (as seperable Hilbert spaces they are isomorphic in a non-canonical way anyway).
In this case, the spectra of $D$ and its closure $\bar D$ in $L^2_{\xi_1}(S)$ and $L^2_{\xi_2}(S)$ are the same. In corollary \ref{cor:HilbertSpaceCompactBaseManifold}, we will see that the space $L^2_\xi(S)$ and the spectrum of $D$ does not depend on $\xi$, if the base manifold $M$ is compact.

After these general results, we turn to explicit calculations of spectra. 
So far, most explicit computations of (subsets of) spectra of the Dirac operator on pseudo-Riemannian manifolds are concerned only with eigenvalues of the Dirac operator. 
For example, the \emph{point spectrum} of the pseudo-Riemannian torus has been computed (\cf \mcitep{166}{Baum1981} for $\specp(\set T^{2,1})$ and \mcitep{95}{Lange} for the general case $\specp(\set T^\pq)$).
In both sources the authors derive that $\bar D$ has empty residual spectrum from an argument involving the symmetry of the point spectrum.
To the best of our knowledge, the only remark in the literature concerned with the continuous spectrum of $\bar D$ can be found in \mcitep{36}{Kunstmann}. There the author states (without proof) the relation $\specc(D)=\set C$ for $\set R^\pq$ and $\specc(D)=\set C\setminus \specp(D)$ for $\set T^\pq$.

In the remaining sections \ref{sec:GeneralizedMultiplicationOperators} -- \ref{sec:SpectrumDiracOperatorTorusProduct} we explicitely compute the spectrum of $D$ for some example manifolds. 
The approach we take is based upon Fourier transform and Fourier series.
Using this method we show that in our example cases the Dirac operator as an unbounded operator in $L^2_\xi(S)$  is unitarily equivalent to a generalized multiplication operator in the Hilbert space $L^2(\Omega,H)$, where $\Omega=\set Z^n$ or $\Omega=\set R^n$ and $H$ is a suitable Hilbert space.

Generalized multiplication operators are introduced in section \ref{sec:GeneralizedMultiplicationOperators}, where we also develop a method to compute their spectra.  
In section \ref{sec:SpectrumDiracOperatorRpq} we compute the spectrum of the Dirac operator of $\set R^\pq$ with indefinite metric of signature $(p,q)$ for the ``canonical'' choice of $\xi:=\linearspan{e_1,\ldots,e_p}$ and the spectrum of the Dirac operator on the flat tori $\set T^\pq$ with indefinite metric of signature $(\pq)$.
Section \ref{sec:SpectrumDiracOperatorTorusProduct} is devoted to the study of the spectrum of the Dirac operator of products of the form $\set T^{1,1}\times F$, where $F$ is a compact even-dimensional Riemannian spin manifold, respectively. All examples share the property that the spectrum of the Dirac operator is $\set C$. On $\set R^\pq$ the spectrum consists only of the continuous spectrum, whereas on the compact manifolds the spectrum consists of continuous and point spectrum.


\section{Notation and prerequisites}
\label{sec:notation}
The purpose of this section is to fix some notation concerning Clifford algebras, spin groups and spinors, spin manifolds and the Dirac operator and its spectrum.
Proofs are found in any textbook on spin geometry, \eg \mcite{Lawson1989}, \mcite{Baum1981} or \mcite{Friedrich2000}.

\subsubsection*{Clifford algebras, spinors and inner products}

Throughout the paper let $p+q=n\geq 2$. 
Let $\Cl_\pq:=\Cl(\set R^n,\sprpqdummy)$ be the Clifford algebra of the bilinear form
			\[
				\sprpq{x,y}:=-x_1y_1 - \ldots - x_py_p + x_{p+1}y_{p+1} + \ldots + x_ny_n \quad \text{for } x,y\in \set R^n,
\]
where we have the relation
\[
	x\cdot y + y\cdot x=-2\sprpq{x,y}\quad \text{for } x,y\in \set R^n.
\]
The usual Euclidean inner product and corresponding norm on $\set R^n$ will be denoted by $\sprRndummy$ and $\normEucldummy$ respectively.
We write $\set R^\pq$ if we endow $\set R^n$ with the bilinear form $\sprpqdummy$. The corresponding pin and spin groups will be denoted by $\Pin(\pq)$ and $\Spin(\pq)$, respectively. The connected component of the neutral element $e\in\Spin(\pq)$ will be denoted by $\Spin_0(\pq)$.

An irreducible representation of the complexified Clifford algebra $\Clc_\pq$ is called a Dirac spinor representation.
Let $e_1,\ldots,e_n$ be the standard basis of $\set R^n$ and let $m=[\f{n}{2}]$. Let furthermore $\Mat{2^m}$ be the space of complex $2^m\times 2^m$-matrices.
The maps
\[
	\begin{alignedat}{2}
		&\Phi_\pq: \Clc_\pq \to \Mat{2^m} & & \quad \text{if } p+q=2m\\
		&\Phi_\pq: \Clc_\pq \to \Mat{2^m} \oplus \Mat{2^m} & & \quad \text{if }  p+q=2m+1 \\
\end{alignedat}
\]
defined by 
\[
	\begin{alignedat}{2}
&\Phi_\pq(e_{j})=\tau(j)\,  E \tensor \ldots \tensor E \tensor U_{\sigma(j)} \tensor \underbrace{T \tensor \ldots \tensor T }_{[\f{j-1}{2}] \text{ times}} & & \quad \text{if } p+q=2m \\
& \Phi_\pq(e_{j})=\begin{cases} (\Phi_{p,q-1}(e_j),\Phi_{p,q-1}(e_j)) \quad \hspace{0.95cm} \text{ if } j=1,\ldots 2m\\ 
	\tau(2m+1)\, (\im \cdot T^ {\tensor m},-\im \cdot T^{\tensor m})	 \quad \text{if } j=2m+1 \\\end{cases} & & \quad \text{if } p+q=2m+1\\
\end{alignedat}
\]
	where	
\[
E=
\left ( \begin{matrix}
1 & 0\\ 
0 & 1\\
\end{matrix}
\right )
\quad
U_1=
\left ( \begin{matrix}
\im & 0\\ 
0 & -\im\\
\end{matrix}
\right )
\quad
U_2=
\left ( \begin{matrix}
0 &\im\\ 
\im & 0\\
\end{matrix}
\right )
\quad
T=
\left ( \begin{matrix}
0 & -\im\\ 
\im & 0\\
\end{matrix}
\right )
\]

and
\[
\tau(j)=\begin{cases} \im \quad j\leq p \\ 1 \quad j>p \end{cases}
\quad 
\sigma(j)=\begin{cases} 1 \quad \text{$j$ odd} \\ 2 \quad \text{$j$ even}\end{cases}
\]

are well-defined isomorphisms of complex algebras. The symbol $\tensor$ denotes the Kronecker product of matrices.
If $n$ is even, the Dirac spinor representation $\hat \kappa_\pq:=\Phi_\pq$ of $\Clc_{pq}$ is defined by the usual representation of $\Mat{2^m}$ on $\Delta_\pq:=\set C^{2^m}$. If $n$ is odd, we set $\hat \kappa^i_\pq:=\text{pr}_i\circ \Phi_\pq$ for $i=1,2$. $\hat \kappa^i_\pq$ again defines a Dirac spinor representation on $\Delta_\pq:=\set C^{2^m}$. Up to isomorphy, these are the only Dirac spinor representations.
If $p+q=2m+1$, we use $\hat \kappa^1_{\pq}$ from now on and let
	\[
		\kappa_\pq := \begin {cases} \hat \kappa_\pq & \text{$p+q$ even} \\ \hat \kappa^1_\pq & \text{$p+q$ odd}.\end{cases}
	\]
Since $\Spin(p,q)\subset \Cl_\pq \subset \Clc_\pq$, by restricting $\hat \kappa_\pq$ ($n$ even) or $\hat \kappa^i_\pq$ ($n$ odd) onto $\Spin(p,q)$, we find representations of the spin group on $\Delta_\pq$, the so called spin representations. 
The restriction of $\kappa_\pq$ onto $\Spin(p,q)$ will be denoted with the same symbol $\kappa_\pq$.
Thus, 
	\[
		\kappa_\pq:\Spin(\pq)\to \GL(\Delta_\pq)
	\]
is a representation of $\Spin(p,q)$ with representation space $\Delta_\pq=\set C^{2^m}$.

From now on let $0<p<n$ with $m=[\f{n}{2}]$ and let $\sprDdummy$ denote the standard positive definite scalar product on $\Delta_\pq=\set C^{2^m}$:
\[
	\sprD{v,w}:=\sum_{i=1}^{2^m} v_i \bar{w_i}\quad \quad\text{for } v,w\in \Delta_\pq.
\]
Let $\tilde K$ denote the maximal-compact subgrup of $\Pin(\pq)$ defined by
\[
	\begin{split}
		\tilde K:=\{y_1\cdot\ldots\cdot y_{k_1}\cdot x_1\cdot \ldots \cdot x_{k_2}|& y_i \in \linearspan{e_1,\ldots,e_p}\cap H^{n-1}_p, \\
		&x_i\in \linearspan{e_{p+1},\ldots,e_{p+q}}\cap S^{n-1}_p, \, k_1,k_2\in \set N\}
\end{split}
\]
where $\cdot$ denotes Clifford multiplication and 
\[
	\begin{split}
		S^{n-1}_p&:=\{x\in \set R^n| \sprpq{x,x}=1\}\\
		H^{n-1}_p&:=\{x\in \set R^n|\sprpq{x,x}=-1\}.\\
	\end{split}
\]
$\tilde K\cap \Spin(\pq)$ is maximal-compact in $\Spin(\pq)$ and $\tilde K_0:= \tilde K \cap \Spin_0(\pq)$ is maximal-compact in $\Spin_0(\pq)$.
$\sprDdummy$ is invariant under the action of $\tilde K$ on $\Delta_\pq$. 
Let $b\in\Clc_\pq$ be defined as
\[
				b:=\begin{cases} e_1\cdot e_2\cdot\ldots\cdots e_p & \quad p=0,1 \text{ mod } 4\\ \im \cdot e_1 \cdot e_2 \cdot \ldots \cdot e_p & \quad p=2,3 \text { mod }4 \end{cases}
\]
The indefinite inner product $\isprdummy$ of signature $(2^{m-1},2^{m-1})$ on the spinor module $\Delta_\pq$ is defined by
\[
	\isprD{v,w}:=\sprD{b\cdot v,w}\quad\text{for } v,w\in\Delta_\pq.
\]
$\isprDdummy$ is invariant under the action of $\Spin_0(\pq)$ on $\Delta_\pq$. 

\subsubsection*{The Dirac operator on pseudo-Riemannian manifolds and its spectrum}

Let $(M,g,Q,\Lambda)$ be a $n$-dimensional time- and space-oriented pseudo-Riemannian 
 spin manifold (\ie $p>1$ and $q>0$) with $Q$ being the spin and $P_M$ being the connected frame bundle of $M$.
For a bundle $(E,M,\pi)$ the set $\Sec{E}$ denotes the smooth sections of the bundle $E$. Let $U\subset M$ be open, then $E_{|U}:=\pi^{-1}(U)$ denotes the restriction of the bundle onto $U$. We use a similar notation for a single point $m\in M$: $E_{|m}:=\pi^{-1}(m)$.
$\Sec{U,E}$ denotes the smooth local sections of the bundle $E_{|U}$.
Finally, for a vector bundle $E$, $\Secc{E}\subset \Sec{E}$ denotes the space of compactly supported sections of $E$.

Let $S:=Q\times_{\kappa_\pq}\Delta_\pq$ be the spinor bundle of $M$.
Let $\nabla^g$ be the Levi-Civita connection on $M$, and let $\nabla^S$ be the induced covariant derivative on the spinor bundle $S$. The operator
	\[
		D: \Sec{S}\xrightarrow{\nabla^S}\Sec{T^*M\tensor S}\isom \Sec{TM \tensor S} \xrightarrow{\mu} \Sec{S}
	\]
	is called Dirac operator of the spinor bundle $S$. The isomorphism between $\Sec{T^*M\tensor S}$ and $\Sec{TM \tensor S}$ is the \textit{musical isomorphism} induced by the metric $g$ and $\mu$ is the Clifford multiplication.


	Let $(P,M,G,\pi)$ be a principal fibre bundle with structure group $G$, let $\rho:G\to \GL(V)$ be a representation of $G$ over $V$, let $\sprdummy_V$ be a $G$-invariant symmetric ($\set K = \set R)$ or hermitian ($\set K = \set C$) scalar product on $V$.
	Then on the associated bundle $E:=P\times_\rho V$ there is a bundle metric given by\footnote{ 
For a principal fibre bundle $P$ over $M$ with structure group $G$ and a representation $\rho:G\to \GL(V)$, an open set $U\subset M$, an associated vector bundle $E:=P\times_\rho V$, a local section $s\in \Sec{U,P}$ and a local section in the vector bundle $X\in \Sec{U,E}$ we introduce the following notation:
\[
	X^s(m):=v(m),\text{ where } X(m)=[(s(m),v(m))] \quadtextf m \in U.
\]
In other words: $X^s:U\to V$ is the representation of the section $X$ \wrt the frame $s$.}
	\[
	g(e,f):=\spr{e^s,f^s}_V\quad \text{for } e,f\in \Sec{U,E}, s \in \Sec{U,P}.
\]
%
Since $M$ ist space- and time-oriented, we can fix a $p$-dimensional time-like subbundle $\xi\subset TM$ and an orthogonal spacelike complement $\eta\subset TM$.
This choice corresponds to a reduction $Q_\xi$ of the $\Spin_0(\pq)$-principal fibre bundle $Q$ to the maximal-compact subgroup $\tilde K_0\subset \Spin_0(\pq)$. We write $P_\xi:=\Lambda(Q_\xi)$. 
The spinor bundle can be represented as 
\[
	S=Q\times_{\kappa_\pq}\Delta_{\pq}=Q_\xi\times_{\kappa_{\pq|\tilde K_0}} \Delta_\pq.
\]
Local frames $e\in \Sec{U,P_\xi}$ or spin frames $s\in \Sec{U,Q_\xi}$ are called \hldef{$\xi$-adapted}. Let
	\[
		\sprB{\varphi,\psi}:=\sprD{\varphi^s,\psi^s}\quad \text{ for a $\xi$-adapted spin-frame }s\in \Sec{U, Q_\xi}
	\]
	and
	\[
		\isprB{\varphi,\psi}:=\isprD{\varphi^s,\psi^s}\quad \text{ for a spin-frame }s\in \Sec{U, Q}.
	\]
We define a Riemannian metric $r_{g,\xi}$ on $TM$ by
	\[
		r_{g,\xi}(X,Y)=\sprRn{X^b,Y^b}\quad \text{ for a $\xi$-adapted local frame $b\in \Sec{U,P_\xi}$}, X,Y \in TM.
	\]
	Furthermore, define the vector bundle homomorphism $J_\xi:S\to S$ locally by Clifford multiplication with the element $b:= i^{\f{p(p-1)}{2}}e_1\cdot\ldots\cdot e_p\in \Clc_\pq$: 
	\[
		(J_\xi \psi)^s:= b \cdot  \psi^s \quadtextf \psi \in S_{|m}, s\in Q_{\xi|m}, m \in M.
	\]
	with $J_\xi^2=\id_S$ and $\isprB{\varphi,\psi}=\sprB{J_\xi\varphi,\psi}=\sprB{\varphi,J_\xi \psi}$.
%
%

	As a pseudo-Riemannian manifold, $(M,g)$ has a canonical volume element $\dvolMg$. By taking the integral over the whole manifold $\sprBdummy$ and $\isprBdummy$ give rise to inner products in the space $\Secc{S}$ of compactly supported spinors as follows:

\[
\begin{split}
\sprL{\cdot,\cdot}: & \Secc{S}\times \Secc{S} \to \set C\\
&\sprL{\psi_1,\psi_2}=\int_M \sprB{\psi_1(x),\psi_2(x)}\,\dvolMg \\
\end{split}
\]

and in the same way
\[
\begin{split}
\isprL{\cdot,\cdot}: & \Secc{S}\times \Secc{S} \to \set C\\
&\isprL{\psi_1,\psi_2}=\int_M \isprB{\psi_1(x),\psi_2(x)}\,\dvolMg. \\
\end{split}
\]

$(\Secc{S},\sprLdummy)$ is a pre-Hilbert space, 
the completion of $\Secc{S}$ with respect to $\normx{\cdot}$ is denoted by $L^2_\xi(S)$.
We use the subscript $\xi$ to emphasize the dependance on the time-like bundle $\xi$.
The indefinite Hermitian inner product $\isprDdummy$ on $\Secc{S}$ is non-degenerate. $J_\xi:\Secc{S}\to\Secc{S}$ as a map between pre-Hilbert spaces is bijective and bounded \wrt $\normx{\cdot}$. It is formally self-adjoint \wrt $\sprLdummy$:
			\[
		\isprL{\varphi,\psi}=\sprL{J_\xi\varphi,\psi}=\sprL{\varphi,J_\xi \psi}\quad \varphi,\psi \in\Secc{S}.
	\]

Using the relation $\isprL{\varphi,\psi}=\sprL{J_\xi \varphi,\psi}$, it is easy to see that $\isprLdummy$ is continuous \wrt the norm topology in $\Secc{S}\times \Secc{S}$. 
Therefore it can be continued to an inner product in $L^2_\xi(S)\times L^2_\xi(S)$. 
$J_\xi$ can be continued to a bijective, continuous self-adjoint operator $J_\xi:L^2_\xi(S)\to L^2_\xi(S)$. 
We conclude that $(L^2_\xi(S),\sprLdummy,\isprLdummy,J_\xi)$ is a Krein space (for definitions and properties of Krein spaces see \mcite{Bognar1974}).

The Dirac operator $D$ with $\dom D:=\Secc{S}\subset L^2_\xi(S)$ is an unbounded operator in the Hilbert space $(L^2_\xi(S),\sprLdummy)$. 
When we want to emphasize the dependence on $\xi$, we write $D_\xi$ or $\bar D_\xi$ when referring to the Hilbert space operators.

The \emph{Riemannian} Dirac operator is symmetric, and if the underlying manifold is complete, it is essentially self-adjoint. For the \emph{pseudo-Riemannian} Dirac operator there are analogous results involving the Krein space structure in $L^2_\xi(S)$, \eg it has been shown that the operator $\im^p D$ is $J$-symmetric.
An operator in a Krein space with fundamental symmetry is called $J-$symmetric ($J$-essentially self-adjoint, $J$-selfadjoint), iff it is symmetric (essentially self-adjoint, self-adjoint) \wrt the indefinite inner product.
We will not further proceed in this direction, although some very interesting results concerning symmetries of the spectrum of the Dirac operator can be deduced from Krein space arguments (\cf \mcite{Bognar1974}, \mcite{Baum1981}, \mcite{Baum1994}).


\section{The dependence of the Dirac operator on the time-like subbundle $\xi$}
\label{sec:QuasiIsometricBundleMetrics}

In this section let $E$ always be a real (or complex) vector bundle over $M$. All bundle metrics on $E$ are Euclidean (Hermitian), \ie in particular they are positive definite.

\begin{definition}
	\label{def:QuasiIsometricBundleMetrics}
	Two Euclidean (Hermitian) bundle-metrics $h_1$ and $h_2$ on $E$ are called \hldef{quasi-isometric} if there is a constant $C\geq 1$ such that for all sections $e\in \Sec{E}$ we have
	\[
		\f{1}{C}h_1(e(x),e(x))\leq h_2(e(x),e(x))\leq C h_1(e(x),e(x))\quad \forall x \in M.
	\]
\end{definition}

A section $e\in \Sec{E}$ is called bounded \wrt a bundle-metric $h$, if the real-valued function $h(e,e)$ is bounded on $M$. It is easy to see that two bundle-metrics $h_1$ and $h_2$ on $E$ are quasi-isometric iff they have the same bounded sections.
As all continuous functions on a compact manifold are bounded, we immediately derive the following proposition: 

\begin{proposition}
\label{prop:CompactBaseManifoldAnyBundleMetricsAreQuasiIsometric}
	If the base manifold $M$ is compact, any two bundle metrics $h_1$ and $h_2$ on $E$ are quasi-isometric.
\end{proposition}

Let $\dvol$  be a volume element on $M$, for example induced by a semi-Riemannian metric on $M$.
For any positive definite Hermitian bundle metric $h$ in $E$ there is a positive definite Hermitian scalar product on $\Secc{E}$ defined by
\[
	\abbildungsdef{\sprdummy_h}{\Secc{E}\times\Secc{E}}{\set C}{(e_1,e_2)_h}{\int_M h(e_1(m),e_2(m))\, \dvol.}
\]
Let $\normdummy_h:=\sqrt{\sprdummy_h}$ be the norm induced on $\Secc{E}$ by $\sprdummy_h$ and let $L^2_h(E)$ be the completion of $\Secc{E}$ \wrt this norm. $(L^2_h(E),\sprdummy_h)$ is a seperable Hilbert space.

\begin{lemma}
	Let the bundle metrics $h_1$ and $h_2$ be quasi-isometric, let $\sprdummy_1$ and $\sprdummy_2$ be the corresponding $L^2$-scalar products on $\Secc{E}$, and let and let $L^2_{h_1}(E)$ and $L^2_{h_2}(E)$ be the completions of $\Secc{E}$ \wrt the respective norms $\norm{\cdot}_1:=\sqrt{\sprdummy_{1}}$ and $\norm{\cdot}_2:=\sqrt{\sprdummy_{1}}$. Then the following assertions hold:
\begin{enumerate}
\item A sequence in $\Secc{E}$ is a Cauchy sequence \wrt $\norm{\cdot}_{1}$ iff it is a Cauchy sequence \wrt $\norm{\cdot}_{2}$.
\item The identity $\Id_{\Secc{E}}: (\Secc{E},\spr{\cdot,\cdot}_{1})\to (\Secc{E},\spr{\cdot,\cdot}_{2})$ considered as a map between pre-Hilbert spaces is a bounded isomorphism with bounded inverse.
\item The map $\Id_{\Secc{E}}: (\Secc{E},\spr{\cdot,\cdot}_{1})\to (\Secc{E},\spr{\cdot,\cdot}_{2})$ from part (b) can be continued to a map $\Id: L^2_{h_1}(E)\to L^2_{h_2}(E)$. It is a well-defined, bounded ismorphism with bounded inverse.
\end{enumerate}
\label{lem:ComparisonHilbertSpaces}
\end{lemma}
\begin{proof}
	Let the constant $C\geq 1$ be chosen as in definition \ref{def:QuasiIsometricBundleMetrics}. 
	In particular, for all $e\in \Secc{E}$ we find 
	\[
		\f{1}{C}\norm{e}_{1}^2=\int_M \f{1}{C} h_1(e(m),e(m))\, \dvol\leq  \int_M  h_2(e(m),e(m))\, \dvol =\norm{e}_{2}^2
\]
and analogously $\norm{e}_{2}^2\leq C \norm{e}_{1}^2$, therefore we conclude $\f{1}{\sqrt{C}} \norm{e}_1\leq \norm{e}_2 \leq \sqrt{C} \norm{e}_1$.
This proves part (a). Then $\Id: (\Secc{E},\spr{\cdot,\cdot}_{1})\to (\Secc{E},\spr{\cdot,\cdot}_{2})$ is bounded by $\sqrt{C}$ and the inverse is bounded by $\f{1}{\sqrt{C}}$ which proves part (b). 
Part (a) also shows that completions of the pre-Hilbert spaces $(\Secc{E},\spr{\cdot,\cdot}_{1})$ and $ (\Secc{E},\spr{\cdot,\cdot}_{2}) $ coincide.
Therefore the identity $\Id_{\Secc{E}}$ and its inverse $\Id^{-1}_{\Secc{E}}$ can be continued to the maps $\Id: L^2_{h_1}(E)\to L^2_{h_2}(E)$ and $\Id^{-1}: L^2_{h_2}(E)\to L^2_{h_1}(E)$. 
Of course $\Id$ and $\Id^{-1}$ are  isomorphisms and again bounded by $\sqrt{C}$ and $\f{1}{\sqrt{C}}$ which proves statement (c).
\end{proof}

\begin{definition}
	Let $G$ be a Lie group and let $\rho:G\to GL(V)$ be a representation of $G$ on the finite-dimensional normed vector space $(V,\normdummy_V)$. We say that the representation of a \hldef{subset} $A\subset G$ on $(V,\normdummy_V)$ is \hldef{bounded by a constant} $C>0$, iff
\[
	\norm{\rho(a)\cdot v}_V\leq C \norm{v}_V \quad \forall v\in V, \forall a \in A.
\]
\end{definition}

The next lemma translates the notion of quasi-isometry into the language of principal fibre bundles.

\begin{lemma}
	\label{lem:QuasiIsometryOfMetricsInPrincipalFibreBundles}
	Let $(P,M,G,\pi)$ be a principal fibre bundle with structure group $G$, let $\rho:G\to \GL(V)$ be a representation of $G$ over $V$ ($\dim V<\infty$), let $\sprdummy_V$ be a positive definite inner product on $V$ and let $E:=P\times_\rho V$ be the associated vector bundle.
	Let $G'\subset G$ be a Lie subgroup such that $\sprdummy_V$ is invariant under $\rho_{|G'}$.
	Then any two reductions $(Q_1,f_1)$ and $(Q_2,f_2)$ of $P$ to $G'$ define two Euclidean (Hermitian) bundle metrics $h_1,h_2$ on $E$ as follows:
	We have
	\[
		Q_1\times_{\rho_{|G'}} V\isom P\times_\rho V \isom Q_2 \times_{\rho_{|G'}} V,
	\]
	and since $G'$ leaves $\sprdummy_V$ invariant, we can define bundle metrics $h_1,h_2$ in the first and third bundle.
Now let
\[
	A_m:=\{x\in G|\text{ there exists } a_1\in Q_{1|m}, a_2\in Q_{2|m} \text{ with } f_1(a_1)\cdot x = f_2(a_2)\}
	\]
	and $A:=\bigcup_{m\in M} A_m$, and let $A^{-1}_m:=\{x^{-1}|x \in A_m\}$ as well as $A^{-1}:=\bigcup_{m\in M} A^{-1}_m$.

	Then $h_1$ and $h_2$ are quasi-isometric iff the representation of both $A\subset G$ and $A^{-1}\subset G$ is bounded on $(V,\sqrt{\sprdummy_V})$. 
\end{lemma}
\begin{proof}
	First let the representations of $A$ and $A^{-1}$ on $(V,\sqrt{\sprdummy_V})$ be bounded by constants $C_1,C_2>0$.
Let $e\in E_{|m}$ and 
$a_1\in  Q_{1|m}, a_2\in Q_{2|m}$, and furthermore $x\in A_m$ with $f_1(a_1)\cdot x=f_2(a_2)$. 
Then $ e^{f_2(a_2)}=x^{-1} \cdot e^{f_1(a_1)}$ and thus
		\[
			h_1(e,e)=\spr{e^{f_1(a_1)},e^{f_1(a_1)}}_V \quad \text{and} \quad h_2(e,e)=\spr{x^{-1}\cdot e^{f_1(a_1)},x^{-1}\cdot e^{f_1(a_1)}}_V
		\]
		as well as
		\[
			h_2(e,e)=\spr{e^{f_2(a_2)},e^{f_2(a_2)}}_V \quad \text{and} \quad h_1(e,e)=\spr{x\cdot e^{f_2(a_2)} ,x\cdot e^{f_2(a_2)}}_V.
		\]
	Since $A_m\subset A$ and $A^{-1}_m\subset A^{-1}$, the representation of $A_m$ and $A^{-1}_m$ is bounded by $C_1$ and $C_2$, respectively. This shows
\[
	h_2(e,e)\leq C_2^2 h_1(e,e) \quad \text{and} \quad h_1(e,e)\leq C_1^2 h_2(e,e).
\]
Thus
\[
	\f{1}{C_2^2} h_2(e,e)\leq h_1(e,e)\leq C^2_1h_2(e,e).
\]
This whole consideration was independent of the point $m\in M$, therefore setting $C:=\max \{C_1^2,C_2^2\}$ we find
\[
	\f{1}{C}h_2(e,e)\leq h_1(e,e)\leq Ch_2(e,e)
\]
for all $e\in \Sec{E}$, and this shows that $h_1$ and $h_2$ are quasi-isometric.

Now let $h_1$ and $h_2$ be quasi-isometric. Suppose the representation of $A$ was not bounded. Then there exist $m_1,m_2,\ldots\in M$ with $x_1\in A_{m_1}, x_2\in A_{m_2},\ldots$ with $\spr{x_j\cdot v_j,x_j\cdot v_j}_V>j^2 \spr{v,v}_V$ for all $j\in \set N$ for $v_1,v_2,\ldots \in V$ with $\spr{v_j,v_j}_V=1$ for all $j\in \set N$.

Choose disjoint open neighbourhoods $U_{j}$ around the $m_j$ and local sections $s_j:U_j\to Q_2$ and smooth maps $\omega_j:M\to[0;1]$ such that $\omega_i(m_j)=1$ and $\supp \omega_{j}\subset U_j$. Let local vector fields $e_j\in \Sec{U_j,E}$ be defined by 
\[
(e_j(m_i))^ {f_2\circ s_j}=v_i.
\]
$e_j\in\Sec{U_j,E}$ is not defined outside of $U_j$, but since $\omega_j\equiv 0$ outside of $U_j$, $\omega_j e_j\in \Sec{E}$ is defined on $M$ itself.
Now let 
\[
	e:=\sum_{j=1}^\infty e_j.
\]
This sum is a well-defined vector field $e\in \Sec{M}$ since for a given $m\in M$ almost all $e_j(m)=0$.
Obviously $h_2(e(m),e(m))\leq 1$ for all $m\in M$, but 
\[
	h_1(e(m_j),e(m_j))=\spr{x_j\cdot v_j,x_j\cdot v_j}_V\geq j^2.
\]
Therefore there cannot exist a constant $C>0$ with $\f{1}{C} h_1(e(m),e(m))\leq h_2(e(m),e(m))$ for all $m\in M$. This shows that $h_1$ and $h_2$ are not quasi-isometric. If $A^{-1}$ is not bounded, the argument works analogoulsy with the roles of $h_1$ and $h_2$ reversed.
\end{proof}

After these preparations, we use the concept of quasi-isometric bundle metrics to derive a geometric condition under which the spectrum of the Dirac operator is the same for two maximal time-like bundles $\xi_1,\xi_2\subset TM$. 
As introduced in section \ref{sec:notation}, let $\sprBdummy[\xi_1],\sprBdummy[\xi_2]$ and $r_{g,\xi_1},r_{g,\xi_2}$ be the positive definite bundle metrics in $S$ and $TM$, respectively.

\begin{theorem}
	\label{thm:DiracSpectrumsCoincide}
	If the two Riemannian metrics $r_{g,\xi_1}$ and $r_{g,\xi_2}$ are quasi-isometric, the following assertions hold:
	\begin{enumerate}	
		\item \label{item:VectorSpacesSame} The $L^2$-spaces $L^2_{\xi_1}(S)$ and $L^2_{\xi_2}(S)$ are \emph{the same} in the following sense: They are the same considered as vector spaces, and the identity $\Id:L^2_{\xi_1}(S)\to L^2_{\xi_2}(S)$ is a bounded isomorphism with bounded inverse.
	 
		\item \label{item:DiagramCommutes} Let $D_{\xi_i}$ and $\bar D_{\xi_i}$ be the Dirac operator and its closure in the Hilbert space $L^2_{\xi_i}(S)$ for $i=1,2$. Then the operators in the Hilbert spaces commute with $\Id$, \ie $(D_{\xi_2} \circ \Id)(\psi)= (\Id \circ D_{\xi_1})(\psi)$ for all $\psi \in \dom{D_{\xi_1}}$ and 
$(\bar D_{\xi_2}\circ \Id)(\psi)= (\Id \circ \bar D_{\xi_1})(\psi)$ for all $\psi \in \dom {\bar D_{\xi_1}}$.
\item The spectrum and all its parts of $D_{\xi_1}$ and $D_{\xi_2}$ as well as  $\bar D_{\xi_1}$ and $\bar D_{\xi_2}$ coincide.
\end{enumerate}	
	
	In particular, the following diagram is commutative:

	
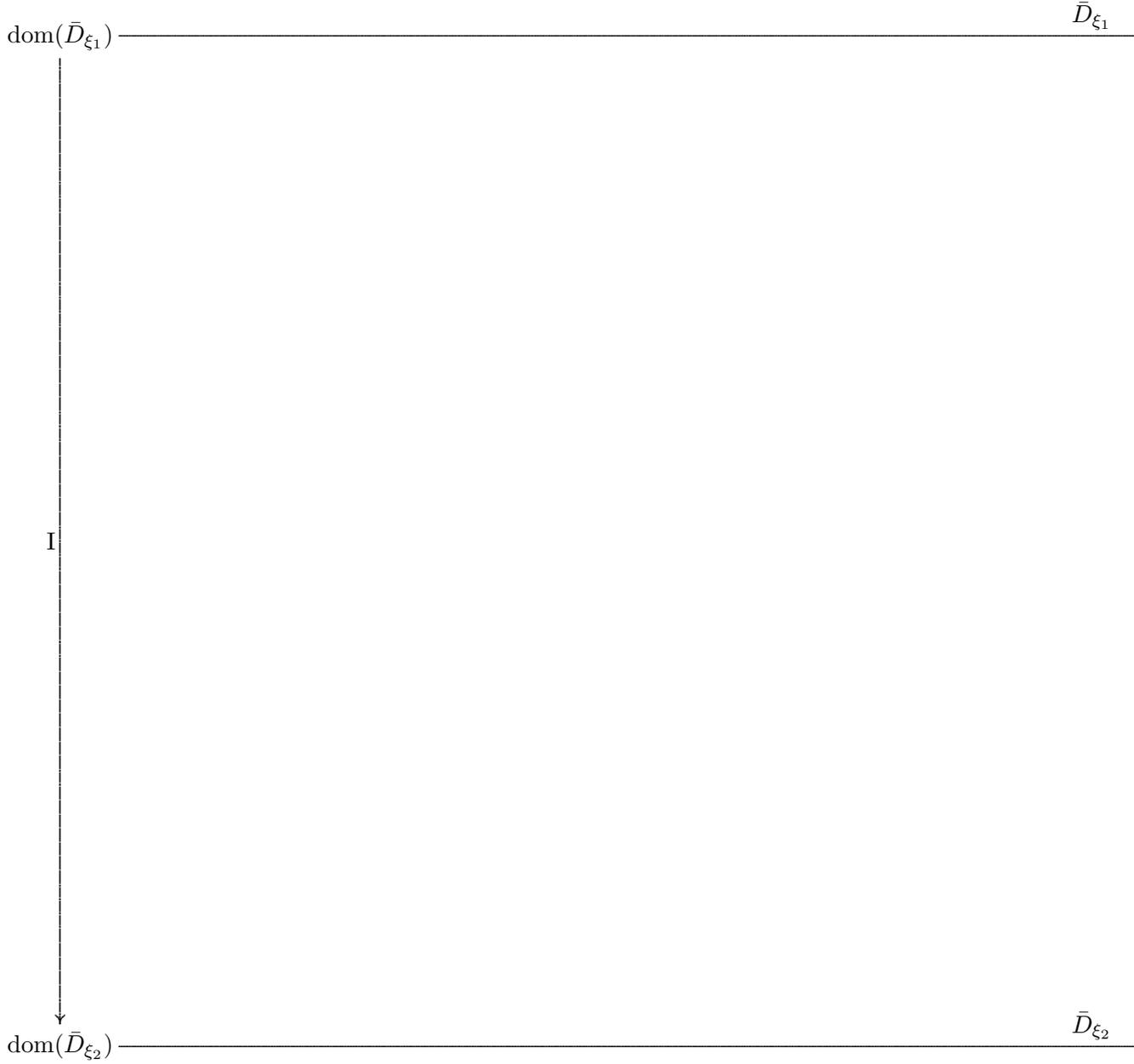
\[
\begindc{\commdiag}
\obj(0,15)[domA1]{$\dom {\bar D_{\xi_1}}$}
\obj(0,0)[domA2]{$\dom {\bar D_{\xi_2}}$}
\obj(30,15)[L21]{$L^2_{\xi_1}(S)$}
\obj(30,0)[L22]{$L^2_{\xi_2}(S)$}
\mor{domA1}{domA2}{$\Id$}[\atright,\solidarrow]
\mor{L21}{L22}{$\Id$}[\atright,\solidarrow]
\mor{domA1}{L21}{$\bar D_{\xi_1}$}[\atleft,\solidarrow]
\mor{domA2}{L22}{$\bar D_{\xi_2}$}[\atleft,\solidarrow]
\enddc
\]	
\end{theorem}

In order to prove the main theorem, we first have to prove the following lemma:

\begin{lemma}
	\label{lem:EquivalenceOfBoundedAction}
	Let $A\subset\Spin(\pq)$. $A$ has bounded representation on $(\Delta_\pq,\sqrt{\sprDdummy})$ iff $\lambda(A)\subset\SO(\pq)$ has bounded representation on $(\set R^n,\normdummy_n)$.
\end{lemma}
\begin{proof}
	To prove this lemma, we will turn both $\Spin(\pq)$ and $\SO(\pq)$ into metric spaces, where the topologies induced by the metrics will coincide with the standard topologies.

\begin{enumerate}
	\item 	To save us some writing, let $M_\pq$ be defined as 
		\[
			M_\pq:=\begin{cases} \Mat{2^m} & p+q=2m \text{ even}\\ \Mat{2^m}\oplus\Mat{2^m} & p+q=2m+1 \text{ odd}.\end{cases}
		\]
		We introduce the norm
		\[
			\normdummy_\pq:=\begin{cases} \normdummy_{\Mat{2^m}} & p+q=2m \text{ even} \\ \max_{i=1,2}(\norm{\pr_i(\cdot)}_{\Mat{2^m}}) & p+q=2m+1 \text{ odd} \end{cases}
		\]
		where $\normdummy_{\Mat{2^m}}$ is the usual operator norm on $\Mat{2^m}$ and $\pr_i:\Mat{2^m}\oplus\Mat{2^m}\to\Mat{2^m}$ denotes the projection on the $i$-th summand.
		
		$(M_\pq,\normdummy_\pq)$ is a normed space, and this gives rise to the metric $d_\pq(A,B):=\norm{B-A}_\pq$ induced by $\normdummy_\pq$. 
		The topology induced by $d_\pq$ is the usual topology of $M_\pq$ as a finite-dimensional complex vector space.
		$\Clc_\pq\supset \Spin(\pq)$ is isomorphic to $M_\pq$ as algebras with $\Phi_\pq:\Clc_\pq\to M_\pq$ being the isomorphism from section \ref{sec:notation}. $\Phi_\pq$ is a homoemorphism \wrt the standard topologies.	
Let $\normdummy_{\Clc_\pq}:=\Phi_\pq^*\normdummy_\pq$ be the pull back of $\normdummy_\pq$ by $\Phi_\pq$. 
%
$(\Clc_\pq,\normdummy_{\Clc_\pq})$ is a normed space. Again, because $\Phi_\pq$ is a homeomorphism, the metric $d_{\Clc_\pq}(a,b):=\norm{b-a}_{\Clc_\pq}$ induces the standard topology on $\Clc_\pq$. 
%
	\item \label{enum:AequivClcBoundedBoundedRep} 
%
	A set $A\subset \Spin(\pq)$ is bounded \wrt the metric $d_{\Clc_\pq}$ iff its representation on $(\Delta_\pq,\sqrt{\sprDdummy})$ bounded by a constant $C>0$.
	This can be seen as follows:
	If $p+q=2m$ is even, the spin representation is given by $\kappa_\pq=\Phi_\pq$, therefore the representation of $A$ on $\Delta_\pq$ is bounded by $C$ iff $\norm{\Phi_\pq(A)}_\pq\leq C$ for all $a\in A$ which in turn is equivalent to the fact that the set $A\subset \Spin(\pq)\subset \Clc_\pq$ is bounded \wrt the metric $d_{\Clc_\pq}$.
	If $p+q=2m+1$ is uneven, things are a bit more complicated:
$\Spin(\pq)$ has two isomorphic Lie group representations on $\Delta_\pq=\set C^{2^m}$, one given by $\hat \kappa^1_\pq=\pr_1\circ \Phi_\pq$, the other given by $\hat \kappa^2_\pq=\pr_2\circ\Phi_\pq$.
	The isomorphism of spin representations  $\Delta_\pq \ni v \mapsto (\pr_2\circ \Phi_\pq)(e_{p+q})v\in \Delta_\pq$ is an isometry \wrt $\sprDdummy$ (see \mcitesatz{1.8}{Baum1981}), therefore $\norm{\pr_1(\Phi_\pq(a))}_{\Mat{2^m}}=\norm{\pr_2(\Phi_\pq(a))}_{\Mat{2^m}}$ for $a\in \Spin(\pq)$.
Thus	
	\[
		\norm{\Phi_\pq(a)}_\pq=\max_{i=1,2}(\norm{\pr_i(\Phi_\pq(a))}_{\Mat{2^m}}) =\norm{\pr_1(\Phi_\pq(a))}_{\Mat{2^m}}
	\]

	Therefore $\norm{\Phi_\pq(a)}_\pq\leq C$ for all $a\in A\subset \Spin(\pq)$, if and only if the representation of $A$ on $\Delta_\pq$ is bounded by $C$. Note again that at this point the set $A$ cannot be an arbitray subset of $\Clc_\pq$ but has to be a subset of $\Spin(\pq)$.

\item	Let $\normdummy_{\mat{n}}$ be the operator norm on the space of real $n\times n$ matrices $\mat{n}$, which again gives rise to a metric $d_{\mat{n}}$ on $\mat{n}$. We denote the restriction onto $\SO(\pq)$ by $d_{\SO(pq)}$.
	The topology induced by $d_{\SO(\pq)}$ is the usual topology on $\SO(\pq)$. 
	
\item \label{enum:AequivBoundedSO} The representation of $B\subset \SO(\pq)$ on $(\set R^n,\normdummy_n)$ is bounded by $C>0$ iff $B$ is bounded \wrt the metric $d_{\SO(pq)}$. This can be seen as follows: The representation of $B$ is bounded by $C>0$ iff
	\[
		\norm{S\cdot v}_n\leq C \norm{v}_n \quad \forall v\in \set R^n, \forall S \in B.
	\]
	This in turn holds true iff $\norm{S}_{\mat{n}}\leq C$ for all $S\in B$, which is equivalent to $B$ being bounded \wrt $d_{\mat{n}}$ and thus also \wrt $d_{\SO{\pq}}$.

\end{enumerate}
	
Now we are able to prove the lemma\footnote{I would like to thank D. Schüth for pointing out a mistake in the original proof.}: Let $A\subset\Spin(\pq)$ have bounded representation on $(\Delta_\pq,\sqrt{\sprDdummy})$. 
By part \ref{enum:AequivClcBoundedBoundedRep} this is the case iff $A$ is bounded \wrt $d_{\Clc_\pq}$. 
	Since $A$ is bounded \wrt $d_{\Clc_\pq}$, there exists a compact set $K\subset\Spin(\pq)$ such that $K\supset A$. Since $K$ is compact and $\lambda$ is continuous, $\lambda(K)$ is compact too and therefore bounded \wrt $d_{\SO(\pq)}$.
	Thus $\lambda(A)$ has bounded representation on $(\set R^n,\normdummy_n)$ by part \ref{enum:AequivBoundedSO}.

	On the other hand, let $\lambda(A)$ have bounded representation on $(\set R^n,\normdummy_n)$. Again by part \ref{enum:AequivBoundedSO} then $\lambda(A)$ is bounded \wrt $d_{\SO(\pq)}$. 
Then there exists a compact set $K\subset \SO(\pq)$ such that $K\supset \lambda (A)$. 
Choose an open cover $(O_\alpha)_{\alpha\in A}$ of $K$ with the property that for every $\alpha\in A$ the set $\lambda^{-1}(O_\alpha)$ consists of exactly two disjoint open sets $U_{\alpha_1}, U_{\alpha_2}\subset \Spin(\pq)$ which are mapped homeomorphically onto $O_\alpha$ by $\lambda$ and the preimage of $\lambda^{-1}(\bar O_a)$ is the disjoint union of $\bar U_{\alpha_1}$ and $\bar U_{\alpha_2}$. This is possible since $\lambda$ is a two-fold covering map, in particular it is a local homeomorphism.
Since $K$ is compact, there exists a finite sub-cover which covers $\lambda(A)$, \ie $\alpha_1,\ldots,\alpha_s\in A$ with $\bigcup_{1\leq i \leq s} O_{\alpha_i}\supset K$.
As closed subsets of the compact set $K$, all the $\bar O_{\alpha_i}$ are compact, and thus the $\bar U_{\alpha_{i_1}},\bar U_{\alpha_{i_2}}$ are compact too.

Obviously the set $\bigcup_{1\leq i \leq s} (\bar U_{\alpha_{i_1}}\cup \bar U_{\alpha_{i_2}})$ is bounded as a finite union of bounded sets, and it covers $\lambda^{-1}(A)$. Thus $\lambda^{-1}(A)$ is bounded \wrt $d_{\Spin(p,q)}$, and using part \ref{enum:AequivClcBoundedBoundedRep} this again shows that $\lambda^{-1}(A)$ has bounded representation on $(\Delta_\pq,\sqrt{\sprDdummy})$.
\end{proof}

We use this lemma to prove the following fundamental proposition:

\begin{proposition}
	\label{prop:EquivalenceOfQuasiIsometry}
	The two Hermitian bundle metrics  $\sprBdummy[\xi_1]$ and $\sprBdummy[\xi_2]$ in the spinor bundle are quasi-isometric iff $r_{g,\xi_1}$ and $r_{g,\xi_2}$ are quasi-isometric.
\end{proposition}

\begin{proof}
	Let $(Q'_{\xi_1},f_1)$ and $(Q'_{\xi_2},f_2)$ be the reductions of the $\Spin_0(\pq)$- bundle $Q$ to the maximal compact subgroup $\tilde K_0\subset \Spin_0(\pq)$ which correspond to the two maximal time-like bundles $\xi_1\subset TM$ and $\xi_2\subset TM$ as explained in section \ref{sec:notation}, and let $\sprBdummy[\xi_1]$ and $\sprBdummy[\xi_2]$ be the corresponding bundle metrics in the spinor bundle $S$.
	
Define for $m\in M$
	\[
		A_m:=\{x \in \Spin_0(\pq)|\text{ there exists } s_1\in Q'_{\xi_1|m}, s_2\in Q'_{\xi_2|m} \text{ with } f_1(s_1)\cdot x = f_2(s_2)\}
	\]

	and $A:=\bigcup_{m\in M} A_m$, and let $A^{-1}_m:=\{x^{-1}|x \in A_m\}$ as well as $A^{-1}:=\bigcup_{m\in M} A^{-1}_m$.

	Let $B_m:=\lambda(A_m)\subset \SO_0(\pq)$, $B^{-1}_m:=\lambda(A^{-1}_m)\subset \SO_0(\pq)$.

	We are now in the situation of lemma \ref{lem:QuasiIsometryOfMetricsInPrincipalFibreBundles}: $\sprBdummy[\xi_1]$ and $\sprBdummy[\xi_2]$ are quasi-isometric iff both $A$ and $A^{-1}$ have bounded representation on $(\Delta_\pq,\sqrt{\sprDdummy})$.
	By lemma \ref{lem:EquivalenceOfBoundedAction} this is equivalent to the fact that that $B$ and $B^{-1}$ have bounded representation on $(\set R^n, \normdummy_n)$.
	By lemma \ref{lem:QuasiIsometryOfMetricsInPrincipalFibreBundles} this is equivalent to the fact that the Riemannian metrics $r_{g,\xi_1}$ and $r_{g,\xi_2}$ are quasi-isometric. 
\end{proof}

We conclude with the proof of the main theorem:

\begin{proof}[proof of theorem \ref{thm:DiracSpectrumsCoincide}]
	Since $r_{g,\xi_1}$ and $r_{g,\xi_2}$ are quasi-isometric, $\sprBdummy[\xi_1]$ and $\sprBdummy[\xi_2]$ are quasi-isometric by proposition \ref{prop:EquivalenceOfQuasiIsometry}. 
	By lemma \ref{lem:ComparisonHilbertSpaces}, $L^2_{\xi_1}(S)$ and $L^2_{\xi_2}(S)$ are isomorphic and $\Id:L^2_{\xi_1}(S)\to L^2_{\xi_2}(S)$ is a bounded isomorphism with bounded inverse, which proves statement \ref{item:VectorSpacesSame}.
	Then \ref{item:DiagramCommutes} is trivial for the Dirac operator, since $\dom {D_{\xi_1}}=\dom{D_{\xi_2}}$ by definition.
	To show \ref{item:DiagramCommutes} for the closure of the Dirac operator, let $x_j$ be a sequence in $\dom{D_{\xi_1}}$ which is converging to $x$ in $ L^2_{\xi_1}(S)$ such that $D_{\xi_1} x_j$ is converging to $y$ in  $ L^2_{\xi_1}(S)$. 
	Then $x_j$ is converging to $x$ in $ L^2_{\xi_2}(S)$ too, and $D_{\xi_2} x_j=D_{\xi_1} x_j$ is converging to $y$ in  $ L^2_{\xi_2}(S)$. This shows \ref{item:DiagramCommutes} and the commutativity of the diagram.
	Part (c) is now a trivial consequence of \ref{item:DiagramCommutes}.
\end{proof}

\begin{corollary}
	\label{cor:HilbertSpaceCompactBaseManifold}
	Let the base manifold $M$ be compact. Then the Hilbert space $L^2_\xi(S)$ as a vector space does not depend on $\xi$. Also, the identity $\Id:L^2_{\xi_1}(S)\to L^2_{\xi_2}$ is a bounded isomorphism with bounded inverse for all choices of $\xi_1,\xi_2$. 
	Furthermore, the spectrum of the Dirac operator $D$ and of its closure $\bar D$ does not depend on $\xi$.
\end{corollary}
\begin{proof}
	If the base manifold is compact, by proposition \ref{prop:CompactBaseManifoldAnyBundleMetricsAreQuasiIsometric} for any maximal time-like subbundles $\xi_1,\xi_2\subset TM$ the two Riemannian metrics $r_{g,\xi_1}$ and $r_{g,\xi_2}$ (and of course also the bundle metrics $\sprBdummy[\xi_1]$ and $\sprBdummy[\xi_2]$ on the spinor bundle) are quasi-isometric, and thus theorem \ref{thm:DiracSpectrumsCoincide} can be applied.
\end{proof}




\section{Generalized multiplication operators and their spectra}
\label{sec:GeneralizedMultiplicationOperators}
In this section generalized multiplication operators are introduced. 
They are operators acting by pointwise multiplication on $L^2$-functions defined on a measure space $\Omega$ with values in a Hilbert space $H$.
\emph{Matrix multiplication operators} are studied in the paper \mcite{Holderrieth1991}, but the author restricts himself to the case where $H$ is finite-dimensional. The infinite-dimensional case is treated in the thesis \mcite{Tomaschewski2003} of S. Tomaschewski, another very condensend overview can be found in the introduction of the article \mcite{Heymann2012}. 

Throughout this section let $(\Omega,\Sigma,\mu)$ be a $\sigma$-finite complete measure space and $(H,\sprdummy_H)$ a seperable Hilbert space. Let $\mathcal L(H)$ denote the space of bounded operators defined on the whole space $H$ and let $\Id:H\to H$ denote the identity.  
For an operator $A$ in $H$ and $\lambda \in \set C$ let $A_\lambda:=\lambda\Id -A$, then the different parts of the spectrum of $A$ are defined as

\begin{tabular}{llp{12cm}}
$\lambda \in \specp(A)$&$\aequiv$ & $\ker{A_\lambda} \neq \{0\} $\\
$\lambda \in \specr(A)$&$\aequiv$ & $\ker{A_\lambda}=\{0\}, \overline{ \ran{A_\lambda}}\neq H $ \\
$\lambda \in \specc(A)$&$\aequiv$ & $\ker{A_\lambda}=\{0\}, \overline{ \ran{A_\lambda}}= H, A_\lambda^{-1}$ is unbounded \\
$\lambda \in \rho(A)$&$ \aequiv$    & $\ker{A_\lambda}=\{0\}, \overline{ \ran{\lambda-A}}= H, A_\lambda^{-1} \text{ is bounded}$\\
\end{tabular}

For a measurable set $A\subset \Omega$ let $\chi_A:\Omega\to\{0,1\}$ denote its characteristic function.
A measurable (complex-valued) function $f$ is called essentially bounded iff it is bounded outside a null set, \ie  there exists a null set $N$ such that $\abs{f_{|\Omega\setminus N}}$ is bounded by $B \in \set R_+$. 
The infimum of all these bounds $B$ is called the \hldef{essential supremum} of $f$.

\begin{definition}
\label{def:GeneralizedMultiplicationOperator}
	Let $M_A:\dom{M_A}\to L^2(\Omega,H)$ be a (possibly unbounded) operator. We call $M_A$ a \hldef{generalized multiplication operator} if 
	\begin{enumerate}
		\item there exists a family $(A(\omega))_{\omega\in\Omega}$ of linear operators with
			\[
				A(\omega):H\to H\quad \text{ if $H$ finite-dimensional}
			\]
			or
			\[
				A(\omega):\dom{A(\omega)}\to H\quad \text{ if $H$ infinite-dimensional}.
			\]
		\item \label{enum:DefiningEquationMultiplicationOperator}
			$M_A$ acts ``by pointwise multiplication with $A(\omega)$'':
			\[
	\begin{split}
	\dom {M_A}=\{f\in L^2(\Omega,H):\, & f(\omega) \in \dom{ A(\omega)} \text{ $\mu$-\ale },\\
& A(\cdot)f(\cdot)) \in L^2(\Omega,H)\},
\end{split}
\]
and
\[
	(M_Af)(\omega)=A(\omega)f(\omega) \quadtextf f\in \dom{M_A}.
\]
	\end{enumerate}
	Conversely, given a family $(A(\omega))_{\omega \in \Omega}$, we define the generalized multiplication operator $M_A:\dom {M_A}\to L^2(\Omega,H)$ associated to $A$ by the relations in \ref{enum:DefiningEquationMultiplicationOperator}.
\end{definition}

\begin{lemma}
	\label{lem:InverseOfMultiplicationOperator}
	Let $M_A$ be a generalized multiplication operator associated to $(A(\omega))_{\omega \in \Omega}$ and let the set
	\[
		\Inv(A):=\{\omega\in \Omega |\, A(\omega) \text{ is invertible} \}
	\]
	be measurable.
	If $\Omega\setminus \Inv(A)$ is a null set, then $M_A$ is invertible and its inverse is the multiplication operator associated to $(A^{-1}(\omega))_{\omega \in \Omega}$. 
\end{lemma}
\begin{proof}
	Let $\Omega\setminus\Inv(A)$ be a null set, and let $M_{A^{-1}}$ be the multiplication operator associated to $(A^{-1}(\omega))_{\omega\in\Omega}$. 
	We first show $\dom{M_{A^{-1}}}=\ran{M_A}$: Let $f\in \dom {M_A}$ and $M_A f \in \ran{M_A}$. Then $(M_Af)(\omega)\in \ran{A(\omega)}$ $\mu$-\ale and therefore $(M_Af)(\omega)\in \dom {A^{-1}(\omega)}$ $\mu$-\ale. 
	Since $A^{-1}(\omega)(M_Af)(\omega)=A^{-1}(\omega)A(\omega)f(\omega)=f(\omega)$ $\mu$-\ale and since $f\in L^2(\Omega,H)$ by definition, this shows $M_Af \in \dom{M_{A^{-1}}}$ and hence $\ran{M_A}\subset \dom{M_{A^{-1}}}$. 

	Using analogous arguments we see that $f':=M_{A^{-1}}f$ is an element of $\dom{M_A}$ for any $f\in \dom{M_{A^{-1}}}$, and $(M_{A}f')(\omega)=f(\omega)$ $\mu$-\ale. Therefore $f\in \ran{M_A}$, $\dom{M_{A^{-1}}}\subset \ran{M_A}$ and thus $\dom{M_{A^{-1}}}= \ran{M_A}$.

	These calculations also show that $M_A \circ M_{A^{-1}} f(\omega)=f(\omega)$ $\mu$-\ale for $f\in \dom {M_{A^{-1}}}$ and $M_{A^{-1}} \circ M_A f(\omega)=f(\omega)$ $\mu$-\ale for $f\in \dom {M_{A}}$.
	Therefore $M_A\circ M_{A^{-1}}\subset \Id$ and $M_{A^{-1}}\circ M_A \subset \Id$ as Hilbert space operators, and finally we conclude $(M_A)^{-1}=M_{A^{-1}}$.
\end{proof}


In the rest of this section we study multiplication operators on the measure spaces $\set N^n, \set Z^n$ and $\set R^n$.

\begin{proposition}
\label{prop:GeneralizedMultiplicationOperatorOnZn}
Let $\Omega=\set Z^n$ or $\Omega=\set N^n$ and let $\mu$ in both cases denote the usual counting measure.
Let $M_A$ be a generalized multiplication operator associated to $(A(k))_{k \in \Omega}$ such that $A(k):\dom{A(k)}\to H$ is densely defined for all $k\in \Omega$.
\begin{enumerate}
	\item $M_A$ is densely defined. \label{enum:ADenselyDefined}
\item \label{enum:MAInvertibleIFFAInvertible} $M_A$ is invertible iff  $A(k)$ is invertible for all $k\in \Omega$.
The inverse $M_A$ is the generalized multiplication operator associated to $A^{-1}(k)$.
\item \label{enum:PointSpectrum} The point spectrum of $M_A$ is given by 
\[
\specp(M_A)=\bigcup_{k \in \Omega} \specp(A(k)).
\]
\item If $\specr(A(k))=\emptyset$ for all $k\in \Omega$, then $\specr(M_A)=\emptyset$. \label{enum:ResidualSpectrum}
\item \label{enum:AdjointOfMultiplicationOperator} The Hilbert space adjoint $(M_A)^*$ of $M_A$  is the  generalized multiplication operator associated to $A^*(k)$.
\item Let $A^*(k)$ be closable for all $k$. Then $\overline{M_A}=M_{\bar A}$. In particular, if $H$ is finite-dimensional or $A(k)$ is closed for all $k\in\Omega$, then $M_A$ is closed itself.
	 \label{enum:MultiplicationOperatorClosed}
 \item $M_A$ is bounded iff $A(k)$ is bounded operator in for all $k\in \Omega$ and there is a constant $B\in \set R_+$ such that $\norm{A(k)}<B$ for all $k\in\Omega$. \label{enum:MABounded}
\end{enumerate}
\end{proposition}
\begin{proof}
	\begin{enumerate}
		\item Let $f\in L^2(\Omega,H)$ and $\Omega_{j}:=\{k \in \Omega: \abs{k_l} \leq j \text{ for all } 1\leq l \leq n \}$.
Then \[
	\mu(\Omega_j)=\begin{cases} (2j+1)^n & \Omega=\set Z^n \\ j^n & \Omega=\set N^n, \end{cases}
\]
$\Omega_{j}\subset \Omega_{j+1}$ and $\Omega=\bigcup_{j=1}^\infty \Omega_{j}$. 
Since $A(k)$ is densely defined for all $k$, we find sequences $(x_j(k))_{j=1,\ldots}\in \dom {A(k)}$ converging to $f(k)$. 
\Wlog we assume $\norm{x_j(k)-f(k)}<\f{1}{(2j+1)^{n}\cdot j}=\f{1}{\mu(\Omega_j)\cdot j}$ ($\Omega=\set Z^n$) or $\norm{x_j(k)-f(k)}<\f{1}{j^{n+1}}=\f{1}{\mu(\Omega_j)\cdot j}$ ($\Omega=\set N^n$).
Let the simple function $f_j\in \dom{M_A}$ be defined by $f_j(k)=\chi_{\Omega_{j}}(k)x_j(k)$. 
Then
\[
	\begin{split} 
		\norm{f-f_j}_{L^2(\Omega,H)}^2&=\int_{\Omega\setminus \Omega_j} \norm{f}^2_H\, d\mu + \sum_{k \in \Omega_j} \norm{f(k)-x_j(k)}^2_H\\
		& < \int_{\Omega\setminus\Omega_{j}} \norm{f}^2_H\, d\mu + \mu(\Omega_j) \cdot \f{1}{\mu(\Omega_j) \cdot j}.
	\end{split}
\]
The second summand tends to $0$ as $j \to \infty$, and the first summand tends to $0$ for $j\to \infty$ because $f\in L^2(\Omega,H)$. This shows that $(f_j)_{j=1,\ldots}\in \dom{M_A}$ converges to $f$, hence $M_A$ is densely defined.

\item One direction is an easy consequence of lemma \ref{lem:InverseOfMultiplicationOperator}: If $A(k)$ is invertible for all $k \in \Omega$, $\Inv(A)=\Omega$ is measurable and $\Omega\setminus\Omega=\emptyset$ is a null set.
	Therefore $M_A$ is invertible with $(M_A)^{-1}=M_{A^{-1}}$. 
	For the other direction, suppose there exists a $k\in \Omega$ such that $A(k)$ is not invertible with $x\in \ker {A(k)}$. Then $\chi_{\{k\}}x$ is a nontrivial element of $\ker {M_A}$, and therefore $M_A$ cannot be invertible.
\item This is a direct consequence of \ref{enum:MAInvertibleIFFAInvertible}: $\lambda-M_A$ is the  maximal generalized multiplication operator associated to $(\lambda-A(k))_{k\in \Omega}$, therefore by \ref{enum:MAInvertibleIFFAInvertible} it has nontrivial kernel iff there exists $k\in \Omega$ such that $\lambda-A(k):\dom{A(k)}\to H$ has nontrivial kernel.
	This is the case iff $\lambda \in \specp(A(k))$. 
\item If $\lambda \in \set C \setminus \specp(M_A)$, then  $\lambda \in \set C \setminus \specp(A(k))$ for all $k\in \Omega$ by \ref{enum:PointSpectrum}.
	Therefore $(\lambda-A)^{-1}(k)$ is defined for all $k\in \Omega$ and densely defined since $\specr(A(k))=\emptyset$ by assumption. Then \ref{enum:MAInvertibleIFFAInvertible} shows that $(\lambda-M_A)^{-1}$ is the multiplication operator associated to $(\lambda-A)^{-1}(k))_{k\in \Omega}$, and therefore $(\lambda-M_A)^{-1}$ is densely defined for all $\lambda \in \set C\setminus \specp(M_A)$ by (a). This shows $\specr(M_A)=\emptyset$.

\item In (a) we showed that $M_A$ is densely defined, thus $(M_A)^*$ is well-defined. Furthermore, $A^*(k)$ is defined for all $k \in \Omega$ since $A(k)$ is densely defined for all $k\in \Omega$. Therefore also $M_{A^*}$ is well-defined.

	First we show $M_{A^*}\subset (M_A)^*$: Let $f\in \dom{M_A}$ and $g\in \dom {M_{A^*}}$. Then
\[
	\begin{split}
	\spr{M_Af,g}_{L^2(\Omega,H)}&=\int_{\Omega} \spr{A(k)f(k),g(k)}_H\, d\mu\\
	&= \int_{\Omega} \spr{f(k),A^*(k) g(k)}_H\, d\mu \\
	& = \spr{f,M_{A*}g}_{L^2(\Omega,H)},
\end{split}
\]
so $g\in \dom {(M_A)^*}$ and hence $M_{A^*}\subset (M_{A})^*$. 
Now let $f\in \dom{(M_A)^*}$. First we have to show that $f(k)\in \dom{A^*(k)}$ for all $k\in \Omega$.
Let $k\in \Omega$, $\chi_{\{k\}}$ be the characteristic function of $\{k\}$ and $x \in \dom{A(k)}$ be arbitrary. Then $\chi_{\{k\}}x\in \dom {M_A}$, and since $f\in \dom{(M_A)^*}$
\[
	\begin{split}
	\spr{A(k)x,f(k)}_H&=\spr{M_A(\chi_{\{k\}}x),f}_{L^2(\Omega,H)}\\
	&=\spr{\chi_{\{k\}}x,(M_A)^*f}_{L^2(\Omega,H)}\\
	&=\spr{x,((M_A)^*f)(k)}_H.
\end{split}
\]
Since $x\in \dom {A(k)}$ was arbitrary, this shows $f(k)\in \dom{A^*(k)}$ with $A^*(k)f(k)=((M_A)^*f)(k)$, hence $\dom{(M_A)^*}\subset \dom {M_{A^*}}$.
The last step is almost trivial now: Let $f\in \dom{(M_A)^*}$. Then for any $\varphi\in \dom{M_A}$ we have
\[
	\begin{split}
	\spr{M_A\varphi,f}_{L^2(\Omega,H)}&=\int_{\Omega} \spr{A(k)\varphi(k),f(k)}_H\, d\mu\\
	&= \int_{\Omega} \spr{\varphi(k),A^*(k)f(k)}_H\, d\mu\\
	&=\spr{\varphi,M_{A^*}f}_{L^2(\Omega,H)}
\end{split}
\]
and therefore $(M_A)^*f=M_{A^*}f$ for all $f\in \dom{(M_A)^*}$ and hence $(M_A)^*=M_{A^*}$.

\item Here we use the following theorem:
	A densely defined operator $B$ in a Hilbert space is closable iff $B^*$ is densely defined, and in this case $\bar B=B^{**}$ (\mcitethm{3.6.3}{Mlak1991}). Since $A(k)$ is closable for all $k$, $\dom {A^*(k)}$ is dense in $H$ for all $k$.
	$(M_A)^*=M_{A^*}$ follows from \ref{enum:AdjointOfMultiplicationOperator}, and this operator is then densely defined by (a). Then $M_A$ in turn is closable and we have
\[
	\bar M_A = (M_A)^{**}=(M_{A^*})^* =M_{A^{**}}.
\]
Now since $A(k)$ is closable for all $k$ we have $\bar A(k)=A^{**}(k)$ for all $k$, and thus $\bar M_A=M_{\bar A}$.
\item Let $M_A$ be bounded and suppose for all $B\in \set R_+$ there exists a $k\in \Omega$ and a $x_k\in \dom{A(k)}$ with $\norm{A(k)x_k}_H>B\norm{x_k}_H$. Then $\chi_{\{k\}}x_k\in \dom {M_A}$ and $\norm{M_A \chi_{\{k\}}x_k}_{L^2(\set Z^n,H}>B \norm{\chi_{\{k\}}x_k}_{L^2(\set Z^n,H)}$, thus $M_A$ cannot be bounded, contradiction.
	For the other direction, assume that $A(k)$ is a bounded operator in $H$ and $\norm{A(k)}<B$ for all $k\in \Omega$.
	Then \[
		\begin{split}
			\norm{M_Af}^2_{L^2(\Omega,H)}&=\int_\Omega \norm{A(k)f(k)}^2_H\, d\mu\\
			&<\int_\Omega B^2 \norm{f(k)}^2_H\, d\mu \\
			&= B^2 \norm{f}^2_{L^2(\Omega,H)}.
		\end{split}
	\]
	for all $f\in \dom {M_A}$. Thus $M_A$ is bounded.
	\end{enumerate}
\end{proof}

\begin{corollary}
	\label{cor:GeneralizedMultiplicationOperatorsOnZn}
Let $\Omega=(\set Z^n,\mu)$ with the usual counting measure.
Let $M_A$ be a generalized multiplication operator associated to $(A(k))_{k \in \set Z^n}$ such that $A(k)$ is densely defined and $\specr(A(k)))=\emptyset$ for all $k\in \set Z^n$. 
Then the spectrum of $M_A$ is
\[
	\sigma(M_A)=\specp(M_A)\cup \specc(M_A)
\]
with 
\[
\specp(M_A)=\bigcup_{k \in \set Z^n} \specp(A(k))
\]
and
\[
	\specc(M_A)=\{\lambda \in \set C: \norm{(\lambda-A(k))^{-1}}\text{ is not bounded on } \set Z^n\}.
\]
Here $\norm{(\lambda-A(k))^{-1}}$ is considered as ``not bounded on $\set Z^n$'' if there exists a $k\in \set Z^n$ such that $(\lambda-A(k))^{-1}:H\to H$ is unbounded already and $\norm{(\lambda-A(k))^{-1}}$ does not even exist for this $k$. 
\end{corollary}

\begin{proof}
	Proposition \ref{prop:GeneralizedMultiplicationOperatorOnZn} \ref{enum:PointSpectrum} and \ref{enum:ResidualSpectrum} imply the first two statements.
	By definition we have
	\[
		\specc(M_A)=\{\lambda \in \set C: (\lambda-M_A)^{-1} \text{ is densely defined and unbounded}\}.
	\]
	From \ref{prop:GeneralizedMultiplicationOperatorOnZn} \ref{enum:MAInvertibleIFFAInvertible} we know that $(\lambda-M_A)^{-1}=M_{(\lambda-A)^{-1}}$. 
	It is clear that this is densely defined, and by \ref{prop:GeneralizedMultiplicationOperatorOnZn} \ref{enum:MABounded} it is bounded iff $(\lambda-A(k))^{-1}(k)$ is a bounded operator for all $k$ and there exists a constant $B\in \set R_+$ with  $\norm{(\lambda-A(k))^{-1}}< B$ for all $k \in \set Z^n$. 
	In this case $\lambda \in \rho(M_A)$, else $\lambda \in \specc(M_A)$.
\end{proof}


The following proposition is the generalization of \mcitethm{III.9.2, p.~134}{Edmunds1987} to our needs. There scalar multiplication operators $(H=\set C$) are treated, but at least the case $\dim H < \infty$ can be treated in the same way. 

\begin{proposition}
	Let $(E,\sprdummy_E)$ be a finite-dimensional Hilbert space. 
	Let $A:\set R^n \to \End(E)$ be measurable and $M_A$ be the multiplication operator associated to $(A(x))_{x\in \set R^n}$. Let $\normdummy_{\End(E)}$ denote the operator norm of $\End(E)$. 
	\begin{enumerate}
		\item $M_A$ is densely defined. \label{enum:MADenselyDefined}
		\item \label{enum:MultiplicationOperatorInverse} The set $\Inv(A)$ from lemma \ref{lem:InverseOfMultiplicationOperator} is measurable. Furthermore $M_A$ is invertible, if $\set R^n\setminus \Inv(A)$ is a null set. The inverse $M_A$ then is a maximal generalized multiplication operator associated to $A^{-1}(x)$.
		\item  \label{enum:MABoundedIFFEssentiallyBounded}
$M_A$ is bounded iff $\norm{A}_{\End(E)}$ is essentially bounded. 
			In this case $\norm{M_A}=\esssup \norm{A}_{\End(E)}$. 		
		\item \label{item:MAclosed} $M_A$ is closed.
	\end{enumerate}
	\label{prop:GeneralizedMultiplicationOperatorOnRn}
\end{proposition}

\begin{proof}
	\begin{enumerate}
		\item Let $f\in L^2(\Omega,E)$ and $\Omega_j:=\{ x \in \set R^n :\norm{A(x)}_{\End(E)}< j\}$ for $j\in \set N$.
	The sets $\Omega_j$ are measurable, since $A$ is measurable and $\normdummy_{\End(E)}$ is continuous.
	Since $A(x)$ as an endomoprhism in the finite-dimensional Hilbert space $E$ is bounded for all $x$, we have $\bigcup_{j=1}^\infty \Omega_j=\set R^n$. 
	Furthermore we have $\Omega_j\subset \Omega_{j+1}$ for all $j$. 
	Let $f_j(x):=f(x) \chi_{\Omega_j}(x)$. Then $f_j\in \dom{M_A}$ because
	\[
		\int_{\set R^n} \norm{A(x) f_j(x)}_E^2 \, dx  \leq j^2 \int_{\Omega_j}  \norm{f(x)}_E^2\, dx< \infty
	\]
	Furthermore, $f_j$ is converging to $f$:
	\[
		\norm{f-f_j}_{L^2(\Omega,E)}^2=\int_{\set R^n \setminus \Omega_j} \norm{f(x)}_E^2\, dx\to 0.
	\]
	This proves that $\dom{M_A}$ is dense.	
		\item $\det(\cdot)$ is a continuous function on $\End(E)$ and obviously $x\in \Inv(A) \aequiv \det(A(x))\neq 0$. Since $A:\set R^n\to \End(E)$ is measurable, the composition $\det \circ A$ is measurable. 
	Therefore $\Inv(A)$ is measurable, and lemma \ref{lem:InverseOfMultiplicationOperator} implies the rest of (b).
	\item 
		Since $A$ is measurable and $\normdummy_{\End(E)}$ is continuous, $\norm{A(\cdot)}_{\End(E)}$ is measurable.

		Assume for some constant $C> 0$ there exists a set $U\subset \set R^n$ such that $\mu(U)>0$ and $\norm{A(x)}_{\End(E)}>C$ for all $x \in U$. \Wlog let $\lambda(U)<\infty$, where $\lambda$ denotes the Lebesgue measure.
		Let $e_1,\ldots,e_j$ be a basis of $E$. For every $x \in U$ we find an $e_l$ such that $\norm{A(x)e_l}> C$ (else $\norm{A(x)}\leq C$).
		For $1\leq l \leq j$ define $U_l:=\{x \in U : \norm{A(x)e_l}>C\}$, then $U=\bigcup_{l=1}^j U_l$.
Since $\lambda(U)>0$ and $U$ is the finite union of the $U_l$, there is at least one $1\leq l \leq j$ such that $0<\lambda(U_l)<\infty$. 
Let $f(x):=\chi_{U_l}(x) e_l$, then $f\in L^2(\set R^n,E)$ and $\norm{f}_{L^2(\set R^n,E)}=\lambda(U_l)>0$. Finally we see that
\[
	\norm{M_Af}_{L^2(\set R^n,E)}^2=\int_{\set R^n} \norm{A(x)f(x)}_E^2\, dx =\int_{U_l} \norm{A(x)e_l}_E^2\, dx>C^2\mu(N_k).
\]
Let $M_A$ now be bounded by $B>0$, and suppose $\norm{A(x)}$  was not essentially bounded by $B$.
Then $\lambda(\{x \in \set R^n: \norm{A(x))}>B\})>0$. But we just saw that in this case $\norm{M_A}>B$, contradiction!
Conversely, if $\norm{A{x}}$ is essentially bounded by $B$, one can directly see that $M_A$ is bounded by $B$, too.	
	
%
%
\item
Since $E$ is finite-dimensional, $A(x)$ is closed (even continuous) for all $x\in \set R^n$. Using \mcitelem{2.3.4}{Tomaschewski2003} this shows that $M_A$ is closed.


	\end{enumerate}

\end{proof}


\section{The spectrum of the Dirac operator on $\set R^{p,q}$ and $\set T^\pq$}
\label{sec:SpectrumDiracOperatorRpq}

In this section we compute the full spectrum of the Dirac operator of $\set R^{\pq}$ in the Hilbert space $L^2_{\xi_e}(S)$ where we choose the maximal time-like subbundle $\xi_e$ to be spanned by the canonical coordinate vector fields $e_1,\ldots,e_p$.
Next, we compute the spectrum of the Dirac operator on the pseudo-Riemannian $n$-dimensional torus $\set T^\pq$. 
Since $\set T^\pq$ is a compact manifold, the spectrum does not depend on the choice of the maximal time-like subbundle $\xi$. 

Both for $\set R^\pq$ and $\set T^\pq$ the computation of the spectrum is a straightforward application of the results of the previous section in spirit of the computation of the spectrum of a constant-coefficient differential operator in $L^2(\set R^n)$ as carried out \eg in \mcitethm{IX.6.2}{Edmunds1987}.

\subsection*{The spectrum of the Dirac operator of $\set R^{p,q}$}

The bundle $Q:=\set R^n\times \Spin_0(\pq)$ defines a spin structure for $\set R^\pq$ with the obvious identifications. 
Since $\set R^n$ is simply-connected, this is the only spin structure.
There is a canonical identification $\Sec{S}\isom C^\infty(\set R^n, \Delta_\pq)$.
With respect to this identification, the Dirac operator has the form
\[
D\psi=\sum_{j=1}^n \kappa_j e_j \cdot \left ( \dd {x^j} \psi \right ), \quad \psi \in C^\infty(\set R^n,\Delta_\pq)
\]
with $\kappa_j:=\sprpq{e_j,e_j}$.
Let a smooth function $A \in C^\infty(\Rpq,\SO_0(\pq))$ be given with $\hat A \in C^\infty(\Rpq,\Spin_0(\pq))$ and $\lambda\circ \hat A = A$.
$A$ induces a maximal time-like bundle $\xi_A$ by
\[
	\xi_A(m):=\linearspan{A(m)e_1,\ldots,A(m)e_p}
\]
			

A quick calculation shows 
\[
	\sprB[\xi_A]{\psi_1,\psi_2}=\sprD{(\hat A(m))^{-1} \psi_1,(\hat A(m))^{-1} \psi_2}.
\]
For the maximal time-like bundle $\xi_e$ spanned by the time-like vectors $e_1,\ldots,e_p$ we find
$\sprB[\xi_e]{\psi_1, \psi_2}=\sprD{ \psi_1, \psi_2}$.
As always, the space of compactly supported smooth spinors $C_0^\infty(\set R^n,\Delta_\pq)$ is a pre-Hilbert space \wrt the inner product
\[
	\sprL[\xi_e]{\psi_1,\psi_2}:=\int_{\set R^n} \sprD{\psi_1(x), \psi_2(x)}\, dx.
\]
and its completion will be denoted by $L^2(\set R^n,\Delta_\pq)$.
Then 
\[
	D:C^\infty_0(\set R^n,\Delta_\pq)\to L^2(\set R^n,\Delta_\pq)
\]
is a densely defined operator in the Hilbert space $L^2(\set R^n,\Delta_\pq)$. 
Let $\FT$ denote the Fourier transform of spinors. $\FT$ is an unitary operator in $L^2(\set R^n,\Delta_\pq)$.
Let $\psi \in C^\infty_0(\set R^n,\Delta_\pq)$ and $x\in \set R^n$, then  
\[
	\begin{split}
		\FT(D\psi)(x)&=\FT(\sum_{j=1}^n \kappa_je_j \cdot \dd {x^j} \psi)(x)\\
		&=\sum_{j=1}^n \kappa_je_j \cdot \FT(\dd {x^j} \psi)(x)\\
		&= \sum_{j=1}^n  \im \kappa_je_j x_j \cdot \FT \psi(x)\\
		&=\im \theta(x) \cdot \FT \psi(x).
	\end{split}
\]
where $\cdot$ denotes the Clifford multiplication and $\theta(x)$ is the reflection of $x$ through the $(n-p)$-plane spanned by the spatial vectors $e_{p+1},\ldots,e_n$.
Let $M_A$ be the generalized multiplication operator associated to 
\[
	\begin{aligned}
		A:&\,\set R^n\to \End(\Delta_\pq)\\
A(x):&\,\Delta_\pq\to \Delta_\pq\\
&\, \psi\mapsto  \im \theta(x) \cdot \psi.\\
\end{aligned}
\]
Note that $A(x)$ depends continuously on $x$ and thus is measurable. 
	Let $M_0$ be the restriction of $M_A$ to 
	\[
		\dom{M_0}=\{ \FT \psi |\, \psi \in C^\infty_0(\set R^n,\Delta_\pq)\}.
	\]
	Then the above calculation shows $D\psi=(\FT^{-1} M_A \FT)\psi$ for $\psi \in C^\infty_0(\set R^n,\Delta_\pq)$ and $D=\FT^{-1}M_0\FT$ as Hilbert space operators. In particular, $D$ and $M_0$ are unitarily equivalent.

\begin{proposition}
	\label{prop:UnitaryEquivalenceOfDirac}
	$M_0$ is closable with $\bar M_0=M_A$. Furthermore, $D$ is closable and for its closure $\bar D$ holds:
	\[
		\begin{aligned}
				\bar D&=\FT^{-1}M_A \FT\\
				\dom{\bar D}&= \{\psi|\, \psi \in L^2(\set R^n, \Delta_\pq) \text{ and } \theta(x)\cdot (\FT\psi)(x) \in L^2(\set R^n,\Delta_\pq)\}.
		\end{aligned}
	\] 
\end{proposition}
\begin{proof}
	Since $M_A$ is closed by proposition \ref{prop:GeneralizedMultiplicationOperatorOnRn} \ref{item:MAclosed}, $M_0$ is closable and $\bar M_0\subset M_A$.
Since $D$ and $M_0$ are unitarily equivalent, also $D$ is closable and $\bar D=\FT^{-1}\bar M_0\FT$. Since $\bar M_0\subset M_A$, it only remains to show that $M_A\subset \bar M_0$. We will prove this in three steps:
\begin{enumerate}[label=(\arabic*)]
\item \label{enum:CompactlySupportedDenseInHM} 
		Let $H(M):=\dom{M_A}$ be the Hilbert space defined by the graph norm
\[
	\norm{x}_{M_A}=\sqrt{\norm{M_A x}_{L^2(\set R^n,\Delta_\pq)}^2+\norm{x}_{L^2(\set R^n,\Delta_\pq)}^2}.
\]
First we show that $C^\infty_0(\set R^n,\Delta_\pq)$ is dense in $H(M)$: 
	Let $\epsilon>0$, and let $\chi_k$ denote the characteristic function of the sets $\Omega_k:=\{x\in\set R^n |\, \norm{x}_n< k\}$.
	For $x\in H(M)$ the sequence $x_k:=\chi_kx$ converges to $x$ in $H(M)$ for $k\to \infty$.
	Choose $K\in \set N$ such that $\norm{x-x_K}_{M_A}<\f{\epsilon}{2}$.
	The function $A:\set R^n\to \End(\Delta_\pq)$ is continuous, therefore $\norm{A(x)}_{\End(\Delta_\pq)}$ is bounded on the compact set $\bar \Omega_K$ by a constant $C_K>0$.
	Since $C^\infty_0(\Omega_K,\Delta_\pq)$ is dense in $L^2(\Omega_K,\Delta_\pq)$, we can choose $x_\epsilon \in C^\infty_0(\Omega_K,\Delta_\pq)$ such that 
	\[
		\sqrt{C_K^2+1}\cdot \norm{x_\epsilon-x_K}_{L^2(\Omega_k,\Delta_\pq)}<\f{\epsilon}{2}.
	\]
	Then we have
	\[
		\begin{split}
		\norm{x-x_\epsilon}_{M_A}&\leq \norm{x-x_K}_{M_A}+\norm{x_K-x_\epsilon}_{M_A}\\
		&<\f{\epsilon}{2}+\left (\int_{\set R^n} \norm{A(y)(x_K-x_\epsilon)(y)}^2_{\Delta_\pq} +\norm{(x_K-x_\epsilon)(y)}^2_{\Delta_\pq}\, dy \right)^{\f{1}{2}}\\
		&=\f{\epsilon}{2}+\left (\int_{\Omega_k} \norm{A(y)(x_K-x_\epsilon)(y)}^2_{\Delta_\pq} +\norm{(x_K-x_\epsilon)(y)}^2_{\Delta_\pq}\, dy \right)^{\f{1}{2}}\\
		&\leq \f{\epsilon}{2}+\left (\int_{\Omega_k} (C_K^2+1) \norm{(x_K-x_\epsilon)(y)}^2_{\Delta_\pq}\, dy \right)^{\f{1}{2}}\\
		&= \f{\epsilon}{2}+\sqrt{C_K^2+1} \norm{x_K-x_\epsilon}_{L^2(\Omega_k,\Delta_\pq)}\\
	&	<\epsilon.
  \end{split}
	\]
Since $x_\epsilon \in C^\infty_0(\Omega_K,\Delta_\pq)\subset C^\infty_0(\set R^n,\Delta_\pq)$ this shows that $C^\infty_0(\set R^n,\Delta_\pq)$ is dense in $H(M)$. 

\item Now we show that $\dom{ M_0}$ is dense in $H(M)$. 
To this end, let $\mathcal S(\set R^n, \Delta_\pq)$ be the Schwartz spinor space of (component-wise) rapidly decreasing spinors and $W^{1,2}(\set R^n,\Delta_\pq)$ be the Sobolev spinor space, again defined component-wise. First we note that $\mathcal S(\set R^n,\Delta_\pq)\subset \dom{M_A}$, since $A(x)$ acts by matrix multiplication with entries polynomial in $x$. 

First we show that $\dom {M_0}$ is dense in $\mathcal S(\set R^n,\Delta_\pq)$ \wrt $\normdummy_{M_A}$: Let $\hat w \in \mathcal S(\set R^n,\Delta_\pq)$ and let $w:=\FT^{-1}\hat w$.
	Since $\FT^{-1}\mathcal S(\set R^n,\Delta_\pq)=\mathcal S(\set R^n,\Delta_\pq)$, we know that $w\in \mathcal S(\set R^n,\Delta_\pq) \subset W^{1,2}(\set R^n,\Delta_\pq)$. 
	The compactly supported spinors $C^\infty_0(\set R^n,\Delta_\pq)$ are dense in $W^{1,2}(\set R^n,\Delta_\pq)$ (\mcitep{222}{Edmunds1987}), thus there exists a sequence $(\phi_k)$ in $C^\infty_0(\set R^n,\Delta_\pq)$ which converges to $w$ in $ W^{1,2}(\set R^n,\Delta_\pq)$. 	
	Therefore $\phi_k\to w$ and $D\phi_k\to D w$ in\footnote{Since $w$ in general is not an element of $\dom D$, since its support need not be compact, by $D w$ we actually mean the differential operator $D$ which of course is defined for the $C^\infty$-function $w$.} $L^2(\set R^n,\Delta_\pq)$, and this in turn means $\FT \phi_k\to \FT w= \hat w$ and $M_A \FT \phi_k \to M_A \FT w=M_A \hat w$ in $L^2(\set R^n,\Delta_\pq)$, which shows $\FT \phi_k \to \hat w$ in $H(M)$.
	Therefore $\dom{M_0}=\FT( C^\infty_0(\set R^n,\Delta_\pq))$ is dense in $\mathcal S(\set R^n,\Delta_\pq)$ \wrt $\normdummy_{M_A}$.

	Since $C^\infty_0(\set R^n,\Delta_\pq) \subset \mathcal S(\set R^n,\Delta_\pq)$, part \ref{enum:CompactlySupportedDenseInHM} shows that $\mathcal S(\set R^n,\Delta_\pq)$ is dense in $H(M)$.
	Therefore $\dom{M_0}$ is dense in $H(M)$, since it is dense in $\mathcal S(\set R^n,\Delta_\pq)$ which in turn is dense in $H(M)$.

\item Since $\dom{M_0}$ is dense in $H(M)$ \wrt $\normdummy_{M_A}$, for all $x\in \dom{M_A}$ there exists a sequence $x_n\in \dom {M_0}$ with $x_n\to x$ and $M_0x_n\to M_A x$ in $L^2(\set R^n,\Delta_\pq)$.
Therefore $x\in \dom {\bar M_0}$  and $M_A x=\bar M_0 x$, which in turn shows $M_A \subset \bar M_0$.
\end{enumerate}
\end{proof}

\begin{lemma}
	\label{lem:InverseOfALambda}
	Let $x\in \set R^n$. The map $A_\lambda(x)\in \End(\Delta_\pq$ defined by $A_\lambda(x):=\lambda-A(x)$ is invertible iff $\lambda^2-\sprpq{x,x}\neq0$, in other words
	\[
	\Inv(A_\lambda)=\{x\in \set R^n|\lambda^2-\sprpq{x,x}\neq 0\}.
\]
For $x\in \Inv(A_\lambda)$ the inverse of $A_\lambda(x)$ is
\[
	(A_\lambda(x))^{-1}=\f{1}{\lambda^2-\sprpq{x,x}}(\lambda+A(x)).
\]
\end{lemma}
\begin{proof}
	Elementary Clifford algebra computations show that 
	\[
		(\lambda-A(x))(\lambda+A(x))=(\lambda+A(x))(\lambda-A(x))=\lambda^2-\sprpq{x,x}
	\]
in the endomorphism algebra $\End(\Delta_\pq)$.
	If $\lambda^2-\sprpq{x,x}\neq0$,
	\[
		(\lambda-A(x))^{-1}(\psi):=\f{1}{\lambda^2-\sprpq{x,x}} (\lambda+A(x)) \psi.
	\]
	is the inverse of $(\lambda-A(x))$. This shows $\Inv(A_\lambda)\supset \{x\in \set R^n| \lambda^2-\sprpq{x,x}\neq 0\}$. 

	To show the converse inclusion, we first prove $\ker{\lambda+A(x)}=\Delta_\pq$ only if $x=0$ and $\lambda=0$.
	Let $\kappa_\pq$ the Dirac spinor representation from section \ref{sec:notation}. Then $\tr(\kappa_\pq(e_j))=0$ for all   $j=1,\ldots,p+q$, since the trace is multiplicative \wrt the Kronecker product and at least one of the matrices $U_1,U_2,T$ from section \ref{sec:notation} with trace $0$ is a factor in the Kronecker product $\kappa_\pq(e_j)$.
	Therefore by linearity also $\tr(\kappa_\pq(x))=0$ for all $x\in \set R^n$.
	Let now $x=\sum_{i_=1}^n x_i e_i$ and $\ker{\lambda +A(x)}= \Delta_\pq$. 
	Then $A(x)=-\lambda \mathbbm{1}$ as matrices, where $\mathbbm{1}$ denotes the identity matrix. 
	Since 
	\[
		0=\tr(\im\kappa_\pq(\theta(x))=\tr(A(x))=\tr(-\lambda\mathbbm{1}) = -\lambda\cdot (\dim \Delta_\pq)
	\]
	we deduce $\lambda=0$ and thus $A(x)=0$. Then by linearity of $\kappa_\pq$ we find $\sum x_i \kappa_\pq(e_i)=0$. Since the $\kappa_\pq(e_i)$ are linearly independent it follows $x=0$.

	Finally we show the inclusion $\Inv(A_\lambda)\subset \{x\in \set R^n|\lambda^2-\sprpq{x,x}\neq 0\}$ by  proving the equivalent statement $\{x\in \set R^n|\lambda^2-\sprpq{x,x}=0\}\subset \set R^n\setminus \Inv(A_\lambda)$.
	Let $x\in \set R^n$ with $\lambda^2-\sprpq{x,x}=0$.	
	If $x=0$, then $\lambda=0$ and obviously $x \not \in \Inv(A_\lambda)$.
	Let now $x\neq 0$. Then $\ker{\lambda+A(x)}\neq \Delta_\pq$ as we have seen above, therefore we can choose $\psi_0\in \Delta_\pq\setminus \ker{\lambda+A(x)}$ and define $\psi_0':= (\lambda + A(x))\psi_0\neq 0$. Then
	\[
		(\lambda-A(x))\psi_0'=(\lambda-A(x))(\lambda+A(x))=(\lambda^2-\sprpq{x,x})\psi_0=0.
	\]

	Therefore $\psi_0'$ is a non-trivial element of $\ker{\lambda-A(x)}$. This shows that $x\not \in \Inv(A_\lambda)$.
\end{proof}

Now we are ready to prove the main theorem:

\begin{theorem}
	\label{thm:SpectrumMinkowski}
	The spectrum of the Dirac operator $D$ of the pseudo-Riemannian manifold $(\set R^n,\sprpqdummy)$ endowed with the canonical spin structure and maximal time-like subbundle $\xi_e$ is the whole complex plane $\set C$, in formulas:
	\[
		\sigma(D)=\sigma(\bar D)=\specc(\bar D)=\set C.
	\]
	In particular, $\specp(\bar D)=\specr(\bar D)=\emptyset$.
\end{theorem}

\begin{proof}
	First note that $\sigma(D)=\sigma(\bar D)$ (\mcitep{105}{Mlak1991}). 
	Since $\bar D$ and $M_A$ are unitarily equivalent by proposition \ref{prop:UnitaryEquivalenceOfDirac}, the spectra and all its parts coincide by \mcitethm{3.4.3}{Mlak1991}. 
	It remains to calculate the spectrum of $M_A$.

	It is easy to see that $M_A$ neither has point nor residual spectrum: Let $\lambda \in \set C$.
	Then $\lambda-M_A=M_{\lambda-A}$ is a multiplication operator as well, and let $A_\lambda$ be defined as in lemma \ref{lem:InverseOfALambda}.
	There, we proved
\[
	\Inv(A_\lambda) = \{x\in \set R^n|\lambda^2-\sprpq{x,x}\neq 0\}.
\]
Then $\set R^n\setminus \Inv(A_\lambda)=\{x\in \set R^n|\lambda^2-\sprpq{x,x}=0\}$, and this is a null set by \mcitelem{IX.6.1}{Edmunds1987}.
Therefore by lemma \ref{lem:InverseOfMultiplicationOperator} the inverse of $\lambda-M_A=M_{\lambda-A}$ is the multiplication operator associated to $A_\lambda^{-1}(x)$. By lemma \ref{lem:InverseOfALambda} we have
\[
	(A_\lambda(x))^{-1}(\psi)= \f{1}{\lambda^2-\sprpq{x,x}} \left ( \lambda \psi + \im \theta(x)\cdot \psi\right ) \quad \text { for } \psi \in \Delta_\pq \text { and } \lambda^2-\sprpq{x,x}\neq 0
\]

Since $(\lambda- M_A)^{-1}=M_{A_\lambda^{-1}}$ is well-defined for all $\lambda \in \set C$, the point spectrum of $M_A$ is empty.
By proposition \ref{prop:GeneralizedMultiplicationOperatorOnRn} \ref{enum:MADenselyDefined}, $\dom{M_{A_\lambda^{-1}}}=\dom{(\lambda-M_A)^{-1}}=\ran{\lambda-M_A}$ is dense. Therefore the residual spectrum $\specr(\bar D)$ is empty, too.
	
It remains to check that $(\lambda-M_A)^{-1}$ is unbounded as an operator in the Hilbert space $L^2(\set R^n,\Delta_\pq)$. 
By proposition \ref{prop:GeneralizedMultiplicationOperatorOnRn} \ref{enum:MABoundedIFFEssentiallyBounded}  it is sufficient to check that $\norm{A_\lambda^{-1}(x)}_{\End(\Delta_\pq)}$ is not essentially bounded.  
This will be done in two steps.
\begin{enumerate}
	\item Since $\normdummy_{\End(\Delta_\pq)}$ is continuous, the function 
		\[
			\Inv(A_\lambda) \ni x \mapsto \norm{A_\lambda^{-1}(x)}_{\End(\Delta_\pq)}
		\]
		is a continuous function defined on an open set.
		We first show that $\norm{A_\lambda^{-1}(x)}_{\End(\Delta_\pq)}$ is unbounded in the usual sense.
		Let first $\lambda=0$, then the continuous function $\norm{A_\lambda(x)}_{\End(\Delta_\pq)}$ is bounded on the compact Euclidean sphere $S^{n-1}$.	
	Since $\f{1}{\sprpq{x,x}}$ is unbounded on $S^{n-1}\cap \{x\in \set R^n| \sprpq{x,x}\neq 0\}$, the function
		\[
			\norm{A_\lambda^{-1}(x)}_{\End(\Delta_\pq)}=\norm{\f{1}{\sprpq{x,x}} A_\lambda(x)}_{\End(\Delta_\pq)}
		\]
		is unbounded on $S^{n-1}\cap \Inv(A_\lambda)$, and therefore $ \norm{A_\lambda^{-1}(x)}_{\End(\Delta_\pq)}$  is unbounded on $\Inv(A_\lambda)$.
		Let now $\lambda \neq 0$ and $X_0:=e_1+e_{p+1}$. 
For $a\in \set R_+$ let $X_a:=aX_0$, then $X_a\in \Inv(A_\lambda)$. 
			By lemma \ref{lem:InverseOfALambda} we have
\[
					 \norm{A_\lambda^{-1}(X_a)}_{\End(\Delta_\pq)}=\norm{\f{1}{\lambda^2}(\lambda +A(X_a)}_{\End(\Delta_\pq)},
			 \]
			and by the homogenity of the norm and the reverse triangle inequality we find the estimate
			 \[\begin{split}
					 &\geq \f{1}{\abs{\lambda^2}}\left | (\norm{A(X_a)}_{\End(\Delta_\pq)}-\abs{\lambda} \right |\\
					 &= \f{1}{\abs{\lambda^2}}\left | ( a\norm{A(X_0)}_{\End(\Delta_\pq)}-|\lambda| \right |.
			 \end{split}
\]
		Since $\norm{A(X_0)}_{\End(\Delta_\pq)}\neq 0$, this is unbounded for $a\to \infty$, therefore $\norm{A_\lambda^{-1}(x)}_{\End(\Delta_\pq)}$ is unbounded on $\Inv(A_\lambda)$.
\item Let now $C\in \set R_+$. Since $\norm{A_\lambda^{-1}(x)}_{\End(\Delta_\pq)}$ is unbounded, we find $x_C\in \Inv{(A_\lambda)}$ with 
			\[
				\norm{A_\lambda^{-1}(x_C)}_{\End(\Delta_\pq)}>C+1.
			\]
			By continuity there exists an open neighbourhood $U_{C}$ of  $x_C$ with non-zero Lebesgue measure and $\norm{A_\lambda^{-1}(x)}_{\End(\Delta_\pq)}>C$ for all $x\in U_{C}$.
			Since $C>0$ was arbitrary, this proves that $\norm{A_\lambda^{-1}(x)}_{\End(\Delta_\pq)}$  cannot be essentially bounded.
\end{enumerate}
\end{proof}

\begin{remark}
Note that $D^2$ is a densely defined differential operator as well, and
\[
	D^2 \psi= - \sum_{j=1}^n \kappa_j \f{\pa^2}{\pa^2 x_j} \psi \quad \text{ for } \psi \in \dom {D^2}:=C^\infty_0(\set R^n,\Delta_\pq).
\]
Therefore $\sigma(\overline{D^2})=\sigma(\bar \Delta)=\set R$, with $\Delta$ being the Laplace-Operator of $\set R^\pq$. Then $\sigma(D)=\set C$, whereas $\sigma(D^2)=\set R$. 
\end{remark}

Now we are considering maximal time-like subbundles other than $\xi_e$, and therefore we also consider different bundle metrics in the spinor bundle.
Let $\xi_A\subset T\set R^\pq$ be another maximal time-like subbundle, and without loss of generality let $\xi_A$ be induced by $A\in C^\infty(\Rpq,\SO_0(\pq))$.
We are interested in geometric conditions on $\xi_A$ under which we can make assertions on the spectrum of $\bar D$ in the space $L^2_{\xi_A}(S)$.

\begin{proposition}
	\label{prop:SpectrumInRn}
	Let a maximal time-like subbundle $\xi_A\subset T\set R^\pq$ be induced by 
\[
	A\in C^\infty(\Rpq,\SO_0(\pq)).
\]
Let the functions $\norm{A(x)}_{\End(\set R^n)}$ and $\norm{A^{-1}(x)}_{\End(\set R^n)}$ be bounded on $\set R^n$.
Let $D_{\xi_A}$ denote the Dirac operator in the Hilbert space $L^2_{\xi_A}(S)$ and $\bar {D}_{\xi_A}$ be its closure. Then they have the same spectral properties as in theorem \ref{thm:SpectrumMinkowski}:	
	\[
		\sigma(D_{\xi_A})=\sigma(\bar D_{\xi_A})=\sigma_c(\bar D_{\xi_A})=\set C.
	\]
\end{proposition}
\begin{proof}
	Since $\norm{A(x)}_{\End(\set R^n)}$ and $\norm{A^{-1}(x)}_{\End(\set R^n)}$ are bounded, the representation of  $\{A(m): m \in \set R^n\}\subset \SO_0(\pq)$ and $\{A^{-1}(m): m \in \set R^n\}\subset \SO_0(\pq)$ is bounded on $(\set R^n,\normdummy_n)$. 
	Lemma \ref{lem:QuasiIsometryOfMetricsInPrincipalFibreBundles} 
	shows that $r_{g,\xi_e}$ and $r_{g,\xi_A}$ are quasi-isometric. 

	Let $D_{\xi_e}$ be the Dirac operator as a operator in $L^2_{\xi_e}(S)=L^2(\set R^n,\Delta_\pq)$. By theorem \ref{thm:DiracSpectrumsCoincide} (c) the spectrum and all its parts of $D_{\xi_e}$ and $D_{\xi_A}$ as well as $\bar D_{\xi_e}$ and $\bar D_{\xi_A}$ coincide. 
\end{proof}

In particular, for any parallel maximal time-like subbundle $\xi'$ on $\Rpq$ the Dirac operator $D_{\xi'}$ has the same spectrum as in theorem \ref{thm:SpectrumMinkowski}.

\begin{exmp}
	\label{example:SpectrumInRn}
	We are now going to discuss one special case in detail: Let $n=2$ and $p=1$. 	
	We have seen that any maximal time-like subbundle $\xi_A$ is induced by a function $A\in C^\infty(\set R^2,\SO_0(1,1)$. Note that in this case the mapping is one-to-one.

Any such $A(m)$ is of the form
\[
	A(m)=\left ( \begin{matrix} \cosh(2a(m)) & \sinh(2a(m)) \\ \sinh(2a(m)) & \cosh(2a(m)) \end{matrix} \right )
\]
with $a\in C^\infty(\set R^2,\set R)$. Therefore we have a one-to-one-correspondence:
\[
	\{\text{maximal time-like subbundles }\xi_A\} \xleftrightarrow{1:1} \{a:a\in C^\infty(\set R^n,\set R)\}.
\]
$a$ can be interpreted as the ``hyperbolic angle'' between $\xi_A$ and the parallel choice $\xi_e$ which of course is induced by $a(x)=0$ for all $x\in \set R^2$.

Now it is easy to see that both $A(m)$ and $A^{-1}(m)$ are bounded if and only if $a(m)$ is bounded. 
For $a(m)\to \infty$, the Euclidean angle between $\xi_A(m):=\linearspan{A(m)e_1}$ and the light ray tends to $0$. Thus we have a very simple geometric interpretation: 
If the \emph{Euclidean} angle between $\xi_A$ and the light cone in $\set R^{1,1}$ is bounded from below (by a constant $>0$), then $D_{\xi_A}$ has the same spectrum as $D_{\xi_e}$.
\end{exmp}

\begin{remark}
	Ideally, one could either prove that the spectrum of the Dirac operator on $\set R^\pq$ is always the whole complex plane $\set C$ -- no matter which maximal time-like subbundle $\xi$ one chooses -- or could give an example of a maximal time-like subbundle $\xi'\subset T\set R^\pq$ such that the spectrum of $D$ in $L^2_{\xi'}(S)$ does not fulfill 
	\[
		\sigma(D_{\xi'})=\sigma(\bar D_{\xi'})=\sigma_c(\bar D_{\xi'})=\set C.
	\]
	Clearly not all $\xi'$ fulfill the conditions of theorem \ref{thm:DiracSpectrumsCoincide}. 
Unfortunately, such a proof or counter-example has yet to be found. 
\end{remark}

\subsection*{The spectrum of the Dirac operator on the torus $\set T^\pq$}
\label{sec:SpectrumDiracOperatorTorus}

Let $\Gamma$ be the discrete subgroup of the group of translations in $\set R^n$ generated by $e_1,e_2,\ldots,e_n$. The $n$-dimensional flat torus $\set T^\pq:=\set R^\pq/(2\pi\Gamma)$ of signature signature $(p,q)$ has $2^n$ spin structures. 
We calculate the spectrum for the trivial spin structure $Q:=\set T^\pq\times \Spin_0(\pq)$ corresponding to $\text{Hom}(\Gamma,\set Z_2)\ni \chi=0$ (\cf \mcitep{113}{Baum1981}).
In this case, $\Sec{S}\isom C^\infty(\set T^\pq, \Delta_\pq)$ and the Dirac operator of a spinor field $\psi\in C^\infty(\set T^\pq,\Delta_\pq)$ has the form
\[
D\psi=\sum_{j=1}^n \kappa_j e_j \cdot \left ( \dd {x^j} \psi \right ).
\]

Since $\set T^\pq$ is compact, the space of square-integrable spinors does not depend on the maximal time-like subbundle $\xi$ and is denoted just by $L^2(S)$ (see corollary \ref{cor:HilbertSpaceCompactBaseManifold}).
We again let $\xi$ be spanned by $e_1,\ldots,e_p$, such that the scalar product for $\psi_1,\psi_2\in L^2(\set T^n,\Delta_\pq)$ is given by 
			\[
				\spr{\psi_1,\psi_2}_{L^2(\set T^n,\Delta_\pq)}:=\f{1}{(2\pi)^n} \int_{\set T^n} \sprD{\psi_1(x),\psi_2(x)}\, dx
			\]
and the Fourier series representation (equally denoted by $\FT$) of a spinor $\psi\in L^2(\set T^n,\Delta_\pq)$
\[
				(\FT\psi)(k)=\f{1}{(2\pi)^n}\int_{\set T^n} e^{-\im\sprRn{k,x}} \psi(x)\, dx\quad \text{for } k\in \set Z^n
\]
	is a unitary transformation between $L^2(\set T^\pq,\Delta_\pq)$ and $L^2(\set Z^n, \Delta_\pq)$.
Let $ \psi \in C^\infty(\set T^n,\Delta_\pq)$ and $k \in \set Z^n$, then
\[	\begin{split}
		\FT(D\psi)(k)&=\FT(\sum_{i=j}^n \kappa_je_j \cdot \dd {x^j} \psi)(k)\\
		& =\sum_{j=1}^n \kappa_je_j \cdot \FT(\dd {x^j} \psi)(k)\\
		&= \sum_{j=1}^n 
		\im  \kappa_je_j k_j \cdot  (\FT \psi)(k)\\
		&=\im  \theta(k) \cdot (\FT \psi)(k).
	\end{split}
\]
Let  $M_A$ be the generalized multiplication operator associated to
\[
	\begin{split}
		A : & \set Z^n \to \End(\Delta_\pq)\\
A(k): &\Delta_\pq \to \Delta_\pq\\
&\psi \mapsto \im \theta(k) \cdot \psi.\\
\end{split}
\]
Let $M_0$ be the restriction of $M_A$ to 
	\[
		\dom{M_0}=\{\FT \psi | \psi \in C^\infty(\set T^n,\Delta_\pq)\}
	\]
	Then the above calculation shows $D\psi=(\FT^{-1} M_A \FT)\psi $ for $\psi \in C^\infty(\set T^\pq,\Delta_\pq)$ and $D=\FT^{-1}M_0\FT$ as Hilbert space operators. In particular, $D$ and $M_0$ are unitarily equivalent.

\begin{proposition}
	\label{prop:UnitaryEquivalenceOfDiracTorus}
$M_0$ is closable with $\bar M_0=M_A$. Furthermore, $D$ is closable and for its closure $\bar D$ holds:
			\[
				\begin{aligned}
					\bar D&=\FT^{-1}M_A \FT\\
				\dom{\bar D}&= \{\psi \in L^2(\set T^\pq, \Delta_\pq)| \theta(k)\cdot \FT( \psi)(k) \in L^2(\set Z^n,\Delta_\pq)\}.
			\end{aligned}
			\] 
\end{proposition}
\begin{proof}
	Since $\Delta_\pq$ is finite-dimensional, $M_A$ is closed by proposition \ref{prop:GeneralizedMultiplicationOperatorOnZn} \ref{enum:MultiplicationOperatorClosed}, and hence $M_0$ is closable and $\bar M_0\subset M_A$.
Because $D=\FT^{-1}M_0\FT$ with $\FT$ being an unitary transformation, $D$ is closable and $\bar D=\FT^{-1}\bar M_0\FT$. It remains to show that $M_A\subset \bar M_0$.
Let $x\in \dom {M_A}$ and let $\chi_j$ denote the characteristic function of the sets $\Omega_j:=\{k\in\set Z^n | \norm{k}_n\leq j\}$. 
Then for $x_j:=\chi_jx$ we have $\FT^{-1}(x_j)\in C^\infty(\set T^n,\Delta_\pq)$ and thus $x_j\in \dom{M_0}$.
Obviously $x_j\to x$, and furthermore $M_0 x_j \to M_A x$ since
	\[
		\norm{M_Ax-M_0x_j}_{L^2}^2=\int_{\set Z^n\setminus \Omega_j} \norm{M_Ax}_{\Delta_\pq}^2 \to 0.
	\]
	This shows $x\in \dom{\bar M_0}$ and $\bar M_0 x=M_Ax$, hence $M_A \subset \bar M_0$.
\end{proof}

\begin{theorem}
	\label{thm:SpectrumTorus}
	The spectrum of the Dirac operator $D$ of $\set T^\pq$  as an unbounded operator in the Hilbert space $L^2(S)$ is the whole complex plane $\set C$, in formulas:
	\[
		\sigma(D)=\sigma(\bar D)=\set C=\specc(\bar D)\cup \specp(\bar D).
	\]
	The point spectrum is
	\[
		\specp(\bar D)=\specp(D)=\{\pm\sqrt{\sprpq{k,k}}|\, k \in \set Z^n\}.
	\]
\end{theorem}
\begin{proof}
	First note that $\sigma(D)=\sigma(\bar D)$ (\mcitep{105}{Mlak1991}). 
	Since $\bar D=\FT^{-1}M_A\FT$ with $\FT$ being an unitary transformation, the spectra and all their parts coincide. 
By corollary \ref{cor:GeneralizedMultiplicationOperatorsOnZn} the spectrum of $M_A$ consists only of the point spectrum $\specp(M_A)$ and the continuous spectrum $\specc(M_A)$:
	\[
		\sigma(M_A)=\specp(M_A)\cup\specc(M_A)	
	\]
with the relations
\[
	\begin{split}
\specp(M_A)&=\bigcup_{k \in \set Z^n} \specp(A(k))\\
	\specc(M_A)&=\{\lambda \in \set C: \norm{(\lambda-A(k))^{-1}}\text{ is not bounded on } \set Z^n\}.\\
\end{split}
\] 
Applying lemma \ref{lem:InverseOfALambda} we find $\lambda \in \specp(A(k))\aequiv \lambda^2-\sprpq{k,k}=0$, therefore
	\[
		\specp(M_A)=\{\pm \sqrt{\sprpq{k,k}}|\, k\in \set Z^n\}.
	\]	
Let $\lambda \in \set C \setminus \specp(M_A)$, in particular $\lambda \neq 0$ from now on.
	Let $a_1:=e_1+e_{p+1}$ and $a_j:=j a_1$ for $j\in \set N$. 
	Choose $\psi_0\in \Delta_\pq$ with $\norm{\theta(a_1)\cdot \psi_0}_{\Delta_\pq}=C>0$  and $\sprD{\psi_0,\psi_0}=1$. By lemma \ref{lem:InverseOfALambda} we have
	\[
		(\lambda-A(a_j))^{-1}\psi_0= \f{\lambda+ \im \theta(a_j)}{\lambda^2} \psi_0.
	\]
	Now using the reverse triangle equation we directly see that
\[
	\begin{split}
		\norm{ (\lambda-A(a_j)^{-1}(\psi_0)}_{\Delta_\pq}&=\norm { \f{\lambda+ \im \theta(a_j)}{\lambda^2} \psi_0}_{\Delta_\pq}\\
		&\geq \left | \norm{\f{\im\theta(a_j)\psi_0}{\lambda^2}}_{\Delta_\pq} - |\f{1}{\lambda}| \right |\\
		&=\left | \f{j C }{\abs{\lambda^2}} - |\f{1}{\lambda}| \right |
\end{split}
\]
The last expression is unbounded for $j\to \infty$, and since $\norm{\psi_0}_{\Delta_\pq}=1$ this means that  $(\lambda-A(a_j))^{-1}$ is unbounded on $\set Z^n$ for all $\lambda \in \set C \setminus \specp(M_A)$. Therefore $\specc(M_A)=\set C \setminus \specp(M_A)$.
\end{proof}

\begin{remark}
	On a compact pseudo-Riemannian spin manifold, given a sequence of spinors $(\psi_j)\in \Secc{S}$ with $D^2\psi_j=0$ and $\norm{D\psi_j}_{L^2(S)}\to \infty$ for $j\to \infty$, the relation $\sigma(\bar D)=\specp(\bar D)\cup \specc(\bar D)=\specap(\bar D)=\set C$ follows from an abstract argument:
	For $\lambda \not \in \specp(\bar D)$ the sequence $h_j:=\f{1}{\norm{\lambda + D\psi_j}_{L^2(S)}}(\lambda + D)\psi_i$ fulfills $\norm{h_j}_{L^2(S)}=1$ for all $j$ and $(\lambda-D)h_j\to 0$ in $L^2(S)$, hence $\lambda$ lies in the approximative spectrum of $\bar D$.
	The previous proof shows that on pseudo-Riemannian spin manifolds such a sequence may exist. One example on $\set T^\pq$ with the trivial spin structure is $\psi_j(x):=e^{-i2\pi\spr{x,a_j}_n}\psi_0$.
	This observation can be used to prove the relation $\sigma(\bar D)=\specp(\bar D)\cup \specc(\bar D)=\set C$ on $\set T^\pq$ for other spin structures than the trivial one by carefully modifying the $\psi_j$. 
\end{remark}
\section{The spectrum of the Dirac operator on a class of compact product manifolds}
\label{sec:SpectrumDiracOperatorTorusProduct}
In this section we compute the spectrum of the Dirac operator of $\set T^{1,1}\times F$ with $F$ being an arbitrary even-dimensional compact Riemannian manifold and $\set T^{1,1}$ being the torus endowed with the pseudo-Riemannian metric of signature $(1,1)$.
This section is technically more complicated than the previous ones, although the main ideas are the same. 
The smooth spinors $\varphi^{\epsilon,k}_l$ defined in lemma \ref{lem:DiracInSplittedBundle} form a Hilbert space basis of $L^2(S)$ and $D\varphi^{\epsilon,k}_l$ is especially easy to compute, which enables us to do Fourier series representation of spinors in $L^2(S)$. 
Along the lines of the previous section we show that $\bar D$ is unitarily equivalent to a generalized multiplication operator. The spectrum then is computed with the methods developed in section \ref{sec:GeneralizedMultiplicationOperators}.

\subsubsection*{Notation and preliminaries: The spin structure on $\set T^{1,1}\times F$}
Let $\set T^{1,1}$ be defined as in the previous section and let $(F,h)$ be a $2N$-dimensional compact Riemannian spin manifold with spin structure $(Q_F,\Lambda_F)$. 
Let $M:=\set T^{1,1}\times F$ be endowed with the product metric $g$ of signature $(1,2N+1)$. 
Let $P_{\set T^{1,1}\times F}$ be the connected frame bundle of $\set T^{1,1}\times F$, 
	$X_1$ and $X_2$ the canonical coordinate vector fields on the product,
$\pr_2:\set T^{1,1}\times F\to F$ the canonical projection and let $P_F$ be the connected frame bundle of $F$.
Consider the reduction of $P_{\set T^{1,1}\times F}$ to the structure group $\SO(2N)$ defined by
\[
	\begin{aligned}
		&\tilde P=\bigcup_{([(x_1,x_2)],f)\in \set T^{1,1}\times F} \tilde P_{|([(x_1,x_2)],f)},\\
		&\tilde P_{|([(x_1,x_2)],f)}:=\{(X_1,X_2,s_1,\ldots,s_{2N})| (s_1,\ldots,s_{2N})=b_m \in P_{F|f}\}
\end{aligned}
\]
where the action of $g\in \SO(2N)$ is explained by $(X_1,X_2,s_1,\ldots,s_{2N})\cdot g := (X_1,X_2,s_1',\ldots,s_{2N}')$ with $(s_1',\ldots,s_{2N}'):=(s_1,\ldots,s_{2N}')\cdot g$.
Let $i_{\SO(2N)}:\SO(N)\to \SO_0(1,2N+1)$ be the embedding defined by
\[
	i_{\SO(2N)}(A)= \left ( \begin{matrix} 1  &  & \\ & 1 & \\ &  & A \\ \end{matrix} \right ).
\]
Then $P \times_{i_{\SO(2N)}} \SO_0(1,2N+1)$ is isomorphic to the connected frame bundle $P_{\set T^{1,1}\times F}$ 
with the isomorphism $\alpha: P\times _{i_{\SO(2N)}}\SO_0(1,2N+1)\to P_{\set T^{1,1}\times F}$ defined by
\[
	[(p,g)]\mapsto g \cdot p
\]
where the action of $g\in \SO_0(1,2N+1)$ on $p\in \tilde P\subset P_{\set T^{1,1}\times F}$ is just the principal fibre action of $\SO_0(1,2N+1)$ on $P_{\set T^{1,1}\times F}$. 
The linear map $i_{\set R^n}:\set R^{2N} \to \set R^{2N+2}$ defined by $i_{\set R^{2N}}(e_j):= e_{j+2}$ for $j=1,\ldots,2N$
induces 
an embedding $i_{\Spin(2N)}:\Spin(2N)\to \Spin(1,2N+1)$. It is easy to see that 
$\lambda \circ i_{\Spin(2N)}=i_{\SO(2N)}\circ \lambda$ and by continuity of all the maps involved the image of $\Spin(2N)$ under $i_{\Spin(2N)}$ is a subset of $\Spin_0(1,2N+1)$. Denote the pullback of $Q_F$ via $\pr_2$ by $\tilde Q:=\pr_2^* Q_F$.

\begin{lemma}
Define the $\Spin_0(1,2N+1)$ principal fibre bundle $Q$ over $\set T^{1,1}\times F$ by $Q:=\tilde Q\times_{i_{\Spin(2N)}}\Spin_0(1,2N+1)$ and let the map $\tilde \Lambda:Q\to \tilde P\times_{i_{\SO(2N)}}\SO_0(1,2N+1)$ be defined by  
\[
	\tilde \Lambda:[(q,g)]\mapsto[(\Lambda_F(q),\lambda(g))].
\]
$\tilde \Lambda$ is well-defined and $(Q,\Lambda)$ with $\Lambda=\alpha\circ \tilde \Lambda$ is a spin structure for $(\set T^{1,1}\times F,g)$.
	\label{lem:SpinStructureOnProduct}
\end{lemma}
\begin{proof}
	$\tilde \Lambda$ is well-defined: Let $h \in \Spin(2N)$, and let $[(q,g)]\in \tilde Q \times_{i_{\Spin(2N)}}\Spin_0(1,2N+1)$. Then 
	\[
	\begin{split}
		\tilde \Lambda ( [(qh,i_{\Spin(2N)}(h^{-1}) \cdot g)] )%
		&=  [(\Lambda_F(qh), \lambda(i_{\Spin(2N)}(h^{-1})\cdot g))]\\
		&= [(\Lambda_F(q)\lambda(h), \lambda(i_{\Spin(2N)}(h^{-1})) \cdot \lambda(g))]\\
		&=[(\Lambda_F(q)\lambda(h), i_{\SO(2N)}(\lambda(h^{-1}) \lambda(g))]\\
		&= [(\Lambda_F (q), \lambda(g))]\\
		&=\tilde \Lambda([(q,g)]).\\
	\end{split}
	\]
It remains to show that $\tilde \Lambda([(q,g)] \cdot h)=\tilde \Lambda([(q,g)])\cdot \lambda(h)$ for $h\in \Spin(1,2N+1)$. 
	Remember that $h$ is acting from the right on the associated bundle:
	\[
	\begin{split}
		\tilde \Lambda([(q,g)]\cdot h)&=\tilde \Lambda([(q,gh)])\\
		&=[( \Lambda_F(q),\lambda(g\cdot h))]\\
		&=[(\Lambda_F(q),\lambda(g))] \cdot \lambda(h)\\
		&=\tilde \Lambda[(q,g)]\cdot \lambda(h)
	\end{split}
	\]
\end{proof}
This construction yields one of the four spin structures on $\set T^{1,1}\times F$.

\subsubsection*{Notation and preliminaries, part II: The spinor bundle and the Dirac operator}
As always let $S:=Q\times_{\kappa_{1,2N+1}} \Delta_{1,2N+1}$ be the spinor bundle of $\set T^{1,1}\times F$. 
First we introduce the vector space isomorphism $\chi:\Delta_{1,2N+1}\to \Delta_{0,2N}\oplus \Delta_{0,2N}$ and show that it is an isomorphism of $\Spin(2N)$ representations. 
Then -- by translating this local isomorphism to a global bundle map -- we prove that the spinor bundle of $\set T^{1,1}\times F$ is isomorphic to $\pr_2^*S_F\oplus \pr_2^*S_F$. This isomorphism then allows us to handle spinor fields in $S$ as tupels of spinor fields in the spinor bundle $S_F$ of $F$ in a very convenient way.

We cite the following lemma from \mciteprop{18}{Baum1997}:
\begin{lemma}
Let $p+q=2m$. For $\epsilon=\pm 1$ let $u(\epsilon)\in \set C^2$ denote the vector 
\[
u(\epsilon)=\f{1}{\sqrt{2}} \left ( \begin{matrix} 1 \\ -\epsilon \im \end{matrix} \right ) \quad \epsilon =\pm1
\]
and $u(\epsilon_1,\ldots,\epsilon_m)=u(\epsilon_1)\tensor\ldots\tensor u(\epsilon_m)$.
Then $\{u(\epsilon_1,\ldots,\epsilon_m)|\prod_{j=1}^m \epsilon_j = \pm 1\}$ is an orthonormal basis of $\Delta_\pq$ \wrt $\sprDdummy$. Furthermore: 
\begin{enumerate}
	\item Every $\psi \in \Delta_{1,2N+1}$ can uniquely be written as $\psi=\psi_1\tensor u(1) + \psi_{-1} \tensor u(-1)$ with $\psi_1,\psi_{-1}\in \Delta_{0,2N}$. 
		The map $\chi: \Delta_{1,2N+1}\to \Delta_{0,2N} \oplus \Delta_{0,2N} $ defined by
\[
	\begin{split}
		\psi \mapsto \chi(\psi):=&(\chi_1(\psi),\chi_{-1}(\psi)):=(\psi_1,\psi_{-1}) \\
		& \text{ for } \psi= \psi_1 \tensor u(1) + \psi_{-1} \tensor u(-1)\\
\end{split}
\]
is a vector space isomorphism.
\item \label{enum:SpinActionCompatible} $\chi$ is compatible with the action of $\Spin(2N)$ on $\Delta_{1,2N+1}$ defined by $\kappa_{1,2N+1}\circ i_{\Spin(2N)}$. Let $s\in \Spin(2N)$ and  $\psi=(\psi_1,\psi_{-1})$, then 
\[
	\chi((\kappa_{1,2N+1}\circ i_{\Spin(2N)})(s) \psi)= (\kappa_{0,2N}(s)\psi_1,\kappa_{0,2N}(s)\psi_{-1}) 
\]
or shorter: $\chi  \circ (\kappa_{1,2N+1}\circ i_{\Spin(2N)})=(\kappa_{0,2N}\oplus \kappa_{0,2N}) \circ \chi$. 
Consequently, the $\Spin(2N)$ representations $\kappa_{1,2N+1}\circ i_{\Spin(2N)}$ and $(\kappa_{0,2N}\oplus \kappa_{0,2N})$ are isomorphic. 

\item \label{enum:PositiveDefiniteScalarProduct}
	For the positive-definite scalare product $\sprD[1,2N+1]{\cdot,\cdot}$ holds:
\[
	\sprD[1,2N+1]{\psi_1,\psi_2}=\sprD[0,2N]{(\chi_1(\psi_1),\chi_1(\psi_2)}+\sprD[0,2N]{\chi_{-1}(\psi_1),\chi_{-1}(\psi_2)}
\]
\item \label{enum:IndefiniteScalarProduct} 
	The indefinite scalar product $\isprdummy_{\Delta_{1,2N+1}}$ is given by
\[
	\isprD[1,2N+1]{\psi_1,\psi_2}=-\sprD[0,2N]{\chi_{-1}(\psi_1),\chi_1(\psi_2)}-\sprD[0,2N]{\chi_{1}(\psi_1),\chi_{-1}(\psi_2)}.
\]
\item \label{enum:CliffordMultiplication} If one considers $\Delta_{0,2N} \oplus \Delta_{0,2N}$ as the representation space of the clifford algebra representation of $\Clc_{1,2N+1}$ via the isomorphism $\chi$, then the clifford multiplication with a vectors from $\set R^{1,2N+1}$ is given by 
\[
	\chi(e_1 \cdot (\psi_1\tensor u(1) + \psi_{-1}\tensor u(-1)))= (-\psi_{-1}, - \psi_1),
\]
\[
	\chi(e_2 \cdot (\psi_1\tensor u(1) + \psi_{-1} \tensor u(-1))  = (-\psi_{-1}, \psi_1),
\]
and for $j>2$
\[
	\chi(e_j \cdot (\psi_1\tensor u(1) + \psi_{-1} \tensor u(-1)))  =  (-e_{j-2}\cdot\psi_1, e_{j-2}\cdot\psi_{-1}).
\]
\end{enumerate}
\label{lem:SpinorRepresentationIsomorphism}
\end{lemma}



Let $\xi$ be the maximal time-like bundle defined by the vector field $X_1$. 
This introduces a positive-definite bundle metric in the spinor bundle $\sprBdummy$.
From now on we always denote a point $m\in \set T^{1,1}\times F$ by $m=([(x_1,x_2)],f)$ with $f\in F, [(x_1,x_2)]\in \set T^{1,1}$.

\begin{lemma}
	Let $S_F:=Q_F\times_{\Spin(2N)} \Delta_{0,2N}$ be the spinor bundle of $F$ and $\pr_2:\set T^{1,1} \times F \to F$ be the canonical projection. 
	\begin{enumerate}
		\item 
			The map
	\[
		B=(B_{1},B_{-1}) :S\to \pr_2^* S_F \oplus \pr_2^* S_F
	\]
	defined for $\epsilon=\pm 1$ by
	\[
		(B_\epsilon(\psi))^q:=\chi_\epsilon(\psi^{\hat q})\quad q\in \Sec{U,\pr_2^*Q_F}, U\subset \set T^{1,1}\times F\text{ open,} 
	\]
	is a well-defined vector bundle isomorphism, where for $q\in \Sec{U,\pr_2^*Q_F}$ we define the spin frame
	\[
		\hat q:=[(q,e)]\in \Sec{U,\pr_2^*Q_F\times_{i_{\Spin(2N)}} \Spin(1,2N+1)}
	\]
\item 
	\[
		\sprB{\psi_1,\psi_2}=\sprBF{B_1(\psi_1),B_1(\psi_2)}+\sprBF{B_{-1}(\psi_1),B_{-1}(\psi_2)},
\]
where $\sprBFdummy$ is the positive-definite bundle metric of the (Riemannian) spinor bundle $S_F$.

	\item \label{enum:IndefiniteSpinorBundle}
		\[
			\isprB{\psi_1,\psi_2}=- \sprBF{B_{-1}(\psi_1),B_1(\psi_2)}-\sprBF{B_1(\psi_1),B_{-1}(\psi_2)}
		\]
	\end{enumerate}	
	
	\label{lem:SpinorBundleIsomorphism}
\end{lemma}
\begin{proof}
	\begin{enumerate}
		\item 
			All we have to show is that the definition of $B$ is independent of the choice of $q$. 
			To this end, let $q'\in \Sec{\pr_2^*Q_F}$, with $q'=q\cdot g$ and $g\in \Spin(2N)$. Then
	\[
		\hat {q'}=[(q',e)]=[(q\cdot g,e)]=[(q,i_{\Spin(2N)}(g))]=[(q,e)]\cdot i_{\Spin(2N)}(g)=\hat q\cdot i_{\Spin(2N)}(g).
	\]
	Therefore $\psi^{\hat {q'}}=i_{\Spin(2N)}(g^{-1}) \psi^{\hat q}$ and finally
	\[
		\chi_\epsilon(\psi^{\hat q'})=\chi_\epsilon(i_{\Spin(2N)}(g^{-1})\psi^{\hat q})=g^{-1}\chi_\epsilon(\psi^{\hat q})
	\]
	because $\chi_\epsilon(i_{\Spin(2N)}(h)\cdot v)=h\cdot \chi_\epsilon (v)$ for all $h \in \Spin(2N)$ by lemma \ref{lem:SpinorRepresentationIsomorphism} \ref{enum:SpinActionCompatible}. This proves that the definition of $B(\psi)$ is in fact independent of the spin frame $q$.
\item Let $q$ and $\hat q$ be defined as in (a). $\hat q$ is $\xi$-adapted, so by definition of $\sprBdummy$ we have
	\[
		\sprB{\psi_1,\psi_2}=\sprD{\psi_1^{\hat q},\psi_2^{\hat q}}\quad \psi_1,\psi_2\in S_{|m}
	\]
	By proposition \ref{lem:SpinorRepresentationIsomorphism} \ref{enum:PositiveDefiniteScalarProduct} we have
	\[
		\begin{split}
		& =\sum_{\epsilon=\pm 1} \sprD[0,2N]{\chi_\epsilon(\psi_1^{\hat q}),\chi_\epsilon(\psi_2^{\hat q})}\\
		&=\sum_{\epsilon=\pm 1} \sprD[0,2N]{(B_\epsilon(\psi_1))^q,(B_\epsilon(\psi_2))^q}\\
		&=\sum_{\epsilon=\pm 1} \sprBF{B_\epsilon(\psi_1),B_\epsilon(\psi_2)}
	\end{split}
\]
\item This follows from lemma \ref{lem:SpinorRepresentationIsomorphism} \ref{enum:IndefiniteScalarProduct}.
\end{enumerate}

\end{proof}

\begin{lemma}
	Let $k\in \set Z^2$ and $\varphi\in \Sec{S_F}$. Let the spinor fields $\varphi^{\epsilon,k}\in \Sec{S}$ be defined by
	\[
		\varphi^{\epsilon,k}([(x_1,x_2)],f):=\begin{cases} B^{-1}( (e^{\im\spr{k,x}_2}\varphi(f),0)) & \epsilon=1\\ B^{-1}( (0,e^{\im\spr{k,x}_2}\varphi(f))) & \epsilon=-1
		\end{cases}
	\]
	
With these definitions the following relations hold for $\epsilon=\pm 1$:
\[
	B_\epsilon(D\varphi^{\epsilon,k})=-\epsilon e^{\im\spr{k,x}} \pr_2^*(D_F \varphi)\quad \text{ and }\quad B_{-\epsilon}(D\varphi^{\epsilon,k})=\im (k_1+\epsilon k_2)e^{\im\spr{k,x}_2}\pr_2^*(\varphi).
	\]
	\label{lem:DiracInSplittedBundle}
\end{lemma}
\begin{proof}
	\begin{enumerate}
		\item 
			Let $U_F\subset F$ be open, $s_F\in\Sec{U_F,Q_F}$, and let $\xi=(s_1,s_2,\ldots,s_{2N})=\Lambda_F \circ s_F$ be the associated local frame in $TF$. Let $s:=\pr_2^*s_F$ be the pulled back frame, \ie $s\in \Sec{\set T^{1,1}\times U_F,\pr_2^*Q}$. Let $\hat s$ be defined as in lemma \ref{lem:SpinorBundleIsomorphism} (a): 
	\[
		\hat s:=[(s,e)]\in \Sec{\set T^{1,1}\times U_F,Q}.
	\]
	With respect to this frame we obviously have
	\[
		(B_\epsilon(\varphi^{\epsilon,k}))^{s}= e^{\im\spr{k,x}_2}(\pr_2^*\varphi)^{s}\quad \text{ and } B_{-\epsilon}(\varphi_{\epsilon,k})=0.
	\]
\item 	The frame in $TM$ associated to $\hat s$ is $\Lambda \circ \hat s=(X_1,X_2,s_1,\ldots,s_{2N})$. A quick calculation using the Koszul formula shows
	\[
		g(\nabla^g_\cdot X_i,\cdot)=0\quad \text{ for }i=1,2
	\]
	\[
		g(\nabla^g_\cdot \cdot ,X_i)=0\quad \text{ for }i=1,2
	\]
and
	\[
		g(\nabla^g_{X_i}\cdot ,\cdot)=0\quad \text{ for }i=1,2
	\]
The only non-vanishing of the local connection forms are
\[
	g(\nabla^g_{s_i}s_j,s_k)=h(\nabla^h_{s_i},s_j,s_k)
\]
Using the local formula of the spinor derivative (\mcitep{216}{Baum1981}) we find 
\[
	(\nabla^S_{X_j}\varphi^{\epsilon,k})^{\hat s}=X_j( (\varphi^{\epsilon,k})^{\hat s} )\quad \text{ for }j=1,2
\]
and for $s_1,\ldots, s_{2N}$:
\begin{equation}
	\label{eqn:DerivativeInCoordinates}
	(\nabla^S_{s_j} \varphi^{\epsilon,k})^{\hat s}=s_j( (\varphi^{\epsilon,k})^{\hat s} ) + \f{1}{2}\sum_{l_1<l_2}^{2N} h(\nabla^h_{s_j}s_{l_1},s_{l_2})e_{l_1+2}\cdot e_{l_2+2}\cdot ( \varphi^{\epsilon,k})^{\hat s}
\end{equation}

\item \label{enum:ProductRule} A quick calculation in coordinates using the product rule shows
\[
	X_j( (\varphi^{\epsilon,k})^{ \hat s})=\im k_j(\varphi^{\epsilon,k})^{\hat s}\quad \text{ for } j=1,2
\]
and therefore 
\[
	\nabla^S_{X_j}\varphi^{\epsilon,k}=\im k_j \varphi^{\epsilon,k}.
\]

\item \label{enum:Derivatives}	Next we want to show
	\[
		B_\epsilon(\nabla^S_{s_j}\varphi^{\epsilon,k})=e^{\im\spr{k,x}_2}\pr_2^*(\nabla^{S_F}_{s_j}\varphi) \quad \text{ and } B_{-\epsilon}(\nabla^S_{s_j}\varphi^{\epsilon,k})=0.
	\]
	By definition we have
	\[
		(B_\epsilon(\nabla^S_{s_j}\varphi^{\epsilon,k}))^{s}=\chi_\epsilon( ( \nabla^S_{s_j}\varphi^{\epsilon,k})^{\hat s} ) 
	\]
	Plugging (\ref{eqn:DerivativeInCoordinates}) into the equation this equals to
	\[
		=\chi_\epsilon(s_j( (\varphi^{\epsilon,k})^{\hat s} ) + \f{1}{2} \sum_{l_1<l_2}^{2N}h(\nabla^h_{s_j} s_{l_1},s_{l_2})\chi_\epsilon(e_{l_1+2}\cdot e_{l_2+2} \cdot (\varphi^{\epsilon,k})^{\hat s} ).
	\]
	Taking together (a), the fact that $s_j(e^{\im\spr{k,x}_2})=0$ and the relations from lemma \ref{lem:SpinorRepresentationIsomorphism} \ref{enum:CliffordMultiplication}, this equals to
	\[
		=e^{\im\spr{k,x}_2} \left ((s_j( (\pr_2^*\varphi)^{s}) ) + \f{1}{2} \sum_{l_1<l_2}^{2N}h(\nabla^h_{s_j} s_{l_1},s_{l_2})e_{l_1}\cdot e_{l_2} \cdot (\pr_2^* \varphi)^{s}) \right ).
	\]
	Note that inside the big brackets nothing depends on $x$, and therefore
	\[
		=e^{\im\spr{k,x}_2} \pr_2^* (\nabla^{S_F}_{s_j}\varphi)^s.
	\]
	An analogous computation for $B_{-\epsilon}$ shows the second statment.
\item Using the local formula of the Dirac operator (\mcitesatz{3.1}{Baum1981}) and \ref{enum:ProductRule}, we have
	\[
		B_\epsilon(D\varphi^{\epsilon,k})=B_\epsilon\left (-\im k_1 X_1 \cdot \varphi^{\epsilon,k}+\im k_2 X_2\cdot \varphi^{\epsilon,k}+\sum^{2N+2}_{j>3} s_{j-2}\cdot \nabla^S_{s_j}\varphi^{\epsilon,k}\right )
	\]
	From the definition of $\varphi^{\epsilon,k}$ and \ref{lem:SpinorRepresentationIsomorphism} \ref{enum:CliffordMultiplication} follows
	\[
		X_1\cdot \varphi^{\epsilon,k}=-\varphi^{-\epsilon,k}\quad \text{and} \quad X_2\cdot \varphi^{\epsilon,k}= \epsilon \varphi^{-\epsilon,k},
	\]
	hence $B_\epsilon(X_1\cdot \varphi^{\epsilon,k}=0$ and $B_\epsilon(X_2\cdot \varphi^{\epsilon,k})=0$.
	Therefore, using \ref{lem:SpinorRepresentationIsomorphism} \ref{enum:CliffordMultiplication}, the whole expression equals to 
\[
	=-\epsilon \sum_{j=1}^{2N} s_j \cdot B_\epsilon( \nabla^S_{s_j}\varphi^{\epsilon,k})%
\]
Using the expressions from part (d) this equals to
\[
	=-\epsilon e^{\im\spr{k,x}_2} \sum_{j=1}^{2N} s_j \cdot \pr_2^*(\nabla^{S_F}_{s_j}\varphi) )
\]
and this in turn equals to
\[
	=-\epsilon e^{\im\spr{k,x}_2} \pr_2^*\left (\sum_{j=1}^{2N} s_j \cdot \nabla_{s_j}^{S_F} \varphi \right )%
	=-\epsilon e^{\im\spr{k,x}_2} \pr_2^*(D_F\varphi).
\]
\item Again using the local formula of the Dirac operator and \ref{enum:ProductRule}, we have
	\[
		B_{-\epsilon}(D\varphi^{\epsilon,k})=B_{-\epsilon} \left (-\im k_1 X_1 \cdot \varphi^{\epsilon,k}+ \im k_2 X_2\cdot \varphi^{\epsilon,k}+\sum^{2N+2}_{j>3} s_{j-2}\cdot \nabla^S_{s_j}\varphi^{\epsilon,k} \right )
	\]
	Since by \ref{enum:Derivatives}
	\[
		B_{-\epsilon}(\nabla^S_{s_j}\varphi^{\epsilon,k})=0\quadtextf j=1,\ldots,2N,
	\]
	we deduce
	\[
			B_{-\epsilon}(D\varphi^{\epsilon,k})=B_{-\epsilon} \left (-\im k_1 X_1 \cdot \varphi^{\epsilon,k} + \im k_2 X_2\cdot \varphi^{\epsilon,k}\right ).
	\]
	Remember $X_1\cdot \varphi^{\epsilon,k}=-\varphi^{-\epsilon,k}$ and $X_2\cdot \varphi^{\epsilon,k}= \epsilon \varphi^{-\epsilon,k}$, then 
this equals to
\[
	\begin{split}
		&=\im k_1B_{-\epsilon}(\varphi^{-\epsilon,k})+\im k_2 B_{-\epsilon} ( \epsilon \varphi^{-\epsilon,k} )\\
		&=\im (k_1+\epsilon k_2) B_{-\epsilon}(\varphi^{-\epsilon,k} )\\
		&=\im (k_1+\epsilon k_2) e^{\im \spr{k,x}_2} \pr_2^*\varphi.
	\end{split}.
\]
	\end{enumerate}

\end{proof}

Since the bundles $\pr_2^*S_F\oplus \pr_2^*S_F$ and $S$ are isomorphic, we can pull back the bundle metric $\sprBdummy$ by $B$. The spaces $\Sec{S}$ and $\Sec{\pr_2^*S_F\oplus \pr_2^* S_F}$ are isomorphic as pre-Hilbert spaces, and therefore $B$ extends to a map to the completions:
\[
	B:L^2(S)\to L^2(\pr_2^*S_F\oplus \pr_2^*S_F)
\]
Of course $B$ and $B^{-1}$ are isometric isomorphisms.

\begin{remark}
	Let us fix some notation: 

	\begin{enumerate}
		\item 
			Let $d_1,d_2\ldots$ be the canonical orthonormal basis of $l^2=L^2(\set N,\set C)$ with $d_l$ being the sequence $d_l:=(\delta_{l,j})_{j=1,\ldots}$, where $\delta_{l,j}$ is the Kronecker delta. Then the set $\{d^\epsilon_l|\, l \in \set N, \epsilon=\pm 1\}$ with
	\[
		d^\epsilon_l:=\begin{cases} (d_l,0) & \epsilon=1 \\ (0,d_l) & \epsilon=-1 \end{cases}
	\]
	forms an orthonormal basis of $l^2\oplus l^2$. 
\item Since $F$ is a compact Riemannian spin manifold, the space $L^2(S_F)$ has a  $L^2$-orthonormal basis of smooth eigenspinors $\varphi_1,\varphi_2,\ldots $ with $D\varphi_l=\lambda_l \varphi_l$. From now on let such a basis be fixed.
\item For $\psi \in \Sec{S}$, let $\psi_\epsilon:=B_\epsilon(\psi)$ be the corresponding sections in $\pr_2^*S_F$.
\end{enumerate}
\end{remark}

\begin{lemma}
	The spaces $L^2(\set T^{1,1},l^2\oplus l^2)$ and $L^2(\pr_2^* S_F \oplus \pr_2^* S_F)$ are unitarily equivalent. The unitary equivalence is given by
	\[
		L:L^2(\pr_2^*S_F\oplus\pr_2^* S_F)\to L^2(\set T^{1,1},l^2\oplus l^2)
	\]
	with $L:= \bar L_0$ being the closure of the operator $L_0$ defined by
	\[
		\dom{L_0}:=\Sec{\pr_2^* S_F\oplus \pr_2^* S_F}
	\]
	and
	\[
		L_0(B(\psi))(x_1,x_2):=\sum_{l=1,\atop \epsilon=\pm 1}^\infty \spr{(\psi_\epsilon(x_1,x_2,\cdot)),\varphi_l}_{L^2(S_F)} \, d^\epsilon_l. 
	\]
	for $B(\psi)=(\psi_1,\psi_{-1})\in\dom{L_0}$ and $[(x_1,x_2)]\in \set T^{1,1}$.
\end{lemma}
\begin{proof}
	In this proof we denote points in $\set T^{1,1}$ with two coordinates $x_1,x_2$ such that $(x_1,x_2)\in \set T^{1,1}$. This is of course meant as an equivalence class \wrt the translative  action of $\Gamma$.

	\begin{enumerate}
		\item We first have to show that $L_0$ is well-defined: Obviously $\psi_\epsilon(x_1,x_2,\cdot)\in \Sec{S_F}\subset L^2(S_F)$ for $B(\psi)\in \dom {L_0}$. Therefore  $\spr{(\psi_\epsilon(x_1,x_2,\cdot)),\varphi_l}_{L^2(S_F)}$ really is a well-defined sequence in $l^2$ and $L_0(B(\psi))(x_1,x_2)$ for $(x_1,x_2)\in \set T^{1,1}$ is a well-defined element of $l^2\oplus l^2$.

		\item $L_0(B(\psi))$ interpreted as a function $\set T^{1,1}\to l^2\oplus l^2$ is measurable:
			By Pettis' theorem it is measurable iff $x'(L_0(B(\psi)))\set T^{1,1}\to \set C$ is measurable for every continuous $x'\in (l^2\oplus l^2)'$. Since $l^2\oplus l^2$ is a Hilbert space, by the Riesz representation theorem we have 
	\[
		x'(\cdot)=\spr{\cdot,\alpha}_{l^2\oplus l^2}
	\]
	with $\alpha \in l^2\oplus l^2$. Without loss of generality assume $\spr{\alpha,d^{-1}_l}_{l^2\oplus l^2}=0$ for all $l$. 
	Then 
\[
	x'(L_0(B(\psi)))=x_\alpha
\]
with 
	\[
		\begin{split}
		x_\alpha:\set T^{1,1} &\to \set C\\
		(x_1,x_2)&\mapsto \spr{\psi_1(x_1,x_2,\cdot),\hat \alpha}_{L^2(S_F)}.
		\end{split}
	\]
	and $\hat \alpha:=\sum_{j=1}^\infty \alpha_j \varphi_j \in L^2(S_F)$. 
Therefore $x'(L_0(B(\psi)))$ is measurable iff $x_\alpha$ is measurable.
\item We show that $x_\alpha$ is measurable by showing that it is continuous:
	Let $\epsilon>0$ and $(x_1,x_2)\in \set T^{1,1}$ be fixed. Obviously the function 
	\[
		\set T^{1,1} \ni (x_1',x_2',f)\mapsto\left ( \psi(x_1',x_2',f)-\psi(x_1,x_2,f) \right )\in S_{F|f}
	\]
	is continuous, and so is $g_{x_1,x_2}(x_1',x_2',f):=\norm{\psi(x_1',x_2',f)-\psi(x_1,x_2,f)}_{S_F|f}$ where $\normdummy_{S_F|f}$ denotes the spinor norm in $S_{F|f}$.
	$g_{x_1,x_2}(x_1,x_2,f)=0$ for all $f\in F$, therefore by the continuity for every $f\in F$ there exists an neighbourhood $U_f\subset \set T^{1,1}\times F$ of $(x_1,x_2,f)$ such that

	\[
		g_{x_1,x_2}(x_1',x_2',f')<\epsilon \quadtextf (x_1',x_2',f')\in U_f.
	\]
	Then for all $f\in F$ the exists a $\delta_f>0$ and an open neighbourhood $U_f'\subset F$ of $f$ such that 
	\[
		(x_1-\delta_f,x_1+\delta_f)\times(x_2-\delta_f,x_2+\delta_f)\times U_f'\subset U_f
	\]
	Since $F$ is compact, already finitely many of the $U_f'$ cover $F$, let therefore $f_1,\ldots,f_n$ be given such that 
	\[
		\bigcup_{j=1}^n U_{f_j}' =F.
	\]
	Let $\delta:=\min_{j=1,\ldots,n} \delta_{f_j}$ and $(x_1',x_2')\in (x_1-\delta,x_1+\delta)\times(x_2-\delta,x_2+\delta)$.
	Then for every $f\in F$ there exists an $f_j$ such that $f\in U_{f_j}'$, and thus we have 
	\[
		(x_1',x_2',f)\in (x_1-\delta,x_1+\delta)\times(x_2-\delta,x_2+\delta)\times U_f'\subset U_f
	\]
	and therefore
	\[
		g_{x_1,x_2}(x_1',x_2',f)<\epsilon
	\]
	for all $f$.
	Therefore
	\[
		\begin{split}
			\abs{x_\alpha((x_1',x_2'))-x_\alpha( (x_1,x_2))}& = \left | \int_{S_F} \spr{\psi_\epsilon(x_1',x_2',f)-\psi_\epsilon(x_1,x_2,f),\hat \alpha (f)}_{S_F|f} \, dh \right |\\
& \leq  \int_{S_F} \left | \spr{\psi_\epsilon(x_1',x_2',f)-\psi_\epsilon(x_1,x_2,f),\hat \alpha (f)}_{S_F|f}\right | \, dh \\
			&\leq \int_{S_F} \norm{\psi_\epsilon(x_1',x_2',f)-\psi_\epsilon(x_1,x_2,f)}_{S_F|f} \cdot \norm{\hat \alpha(f)}_{S_F|f}\, dh \\
			&=\int_{S_F} g_{x_1,x_2}(x_1',x_2',f)\cdot \norm{\hat \alpha(f)}_{S_{F|f}}\, dh \\
			&< \epsilon \int_{S_F} \norm{\hat \alpha(f)}_{S_{F|f}}\, dh\\
			&=\epsilon C\\
		\end{split}
		\]
		for $ (x_1',x_2') \in (x_1-\delta,x_1+\delta)\times(x_2-\delta,x_2+\delta)$ with $C:=\int_{S_F} \norm{\hat \alpha(f)}_{S_{F|f}}\, dh$.
		Thus $x_\alpha$ is continuous in $(x_1,x_2)$. Since this consideration works for all $(x_1,x_2)\in \set T^{1,1}$, $x_\alpha$ is continuous and hence measurable. 
	\item Since this holds true for all $\hat \alpha \in L^2(S_F)$, this proves that $L_0(B(\psi))$ is a measurable function. 
		To show that $L_0$ is well-defined, it remains to show that the $L_0(B(\psi))\in L^2(\set T^{1,1},l^2\oplus l^2)$.
	It is easy to see that
\[
	\begin{split}
	\norm{L_0(B(\psi))(x_1,x_2),L_0(B(\psi))(x_1,x_2)}_{l^2\oplus l^2}&=\sum_{l=1,\atop \epsilon=\pm 1}^\infty \left | \spr{\psi_\epsilon(x_1,x_2,\cdot),\varphi_l}_{L^2(S_F)} \right |^2\\
	&=\sum_{\epsilon=\pm 1} \norm{\psi_\epsilon(x_1,x_2,\cdot)}_{L^2(S_F)}^2.
\end{split}
\]
and therefore
\[
	\begin{split}
		\norm{L_0(B(\psi))}^2_{L^2(\set T^{1,1},l^2\oplus l^2)}=&\int_{\set T^{1,1}} \norm{L_0(B(\psi))(x_1,x_2),L_0(B(\psi))(x_1,x_2)}_{l^2\oplus l^2}^2\, dx_1dx_2\\
		=&\int_{\set T^{1,1}} \sum_{\epsilon=\pm 1}\norm{\psi_\epsilon(x_1,x_2,\cdot)}^2_{L^2(S_F)}\, dx_1dx_2\\
		=&\int_{\set T^{1,1}\times F} \sum_{\epsilon=\pm 1} \sprBF{\psi_\epsilon(x_1,x_2,f),\psi_\epsilon(x_1,x_2,f)}\, dx_1dx_2 dh\\
		=&\int_{\set T^{1,1}\times F} \sprB{\psi(x_1,x_2,f),\psi(x_1,x_2,f)}\, dx_1dx_2 dh\\
		=&\norm{\psi}_{L^2_\xi(S)}<\infty
	\end{split}
\]
where the last equality follows from lemma \ref{lem:SpinorBundleIsomorphism} (b). This both shows that $L_0(B(\psi))\in L^2(\set T^{1,1},l^2\oplus l^2)$ and that $L_0$ is an isometry (and therefore bounded). Now $L_0$ is bounded (and thus closable) and obviously densely defined, therefore its closure $L:=\bar L_0$ is a bounded operator defined on the whole space.

\item Of course $L$ is an isometry as well, and  therefore it is injective. The surjectivity of $L$ can be seen as follows: For $f\in L^2(\set T^{1,1}, l^2\oplus l^2)$ Pettis' theorem shows that the $f^\epsilon_l(x_1,x_2):=\spr{f(x_1,x_2),d^\epsilon_l}_{l^2\oplus l^2}$ are Lebesgue-measurable and  $f^\epsilon_l\in L^2(\set T^{1,1})$. 
		This immediately proves that for all $l$ the function $f^\epsilon_l(x_1,x_2)\varphi_l\in L^2(\pr_2^* S_F)$. Therefore the \emph{finite} sums
			\[
				s_j:=\left (\sum^j_{l=1\atop } f^1_l(x_1,x_2)\varphi_l,\sum^j_{l=1} f^{-1}_l(x_1,x_2)\varphi_l \right ) 
			\]
			lie in  $L^2(\pr_2^* S_F\oplus \pr_2^* S_F)$ as well, and since the $s_j$ form form a Cauchy sequence, they converge to $\hat f$ in the complete space $L^2(\pr_2^* S_F\oplus pr_2^*S_F)$. Since $L(\hat f)=f$, this shows that $L$ is surjective.
	\end{enumerate}
\end{proof}

\begin{remark}
	\label{rem:FourierSeriesOfSpinorsHilbertSpace}
	Fourier coefficients of functions $f\in L^2(\set T^{1,1},H)$ with $(H,\sprdummy_H)$ being a seperable Hilbert space are defined as in \mcite[chapter 5]{Nagy2010}, \ie for any ortonormal Hilbert space basis $b_1,b_2,\ldots$ let and $e^k_l:=e^{\im\spr{k,x}_2}b_l\in L^2(\set T^{1,1},H)$. 
	Then the set $B:=\{e^k_l| k \in \set Z^2, l\in \set N\}$ is a Hilbert space orthonormal basis of $L^2(\set T^{1,1},H)$. 
	The map $\FT$ defined by
	\[\begin{split}
			\FT:B&\to L^2(\set Z^2,H)\\
			e^k_l&\mapsto \chi_{k}b_l\\
	\end{split}
	\]
	can be extended to a unitary transformation $\FT:L^2(\set T^{1,1},H)\to L^2(\set Z^2,H)$, 
	which does not depend on the basis $B$.
\end{remark}


\begin{lemma}
	\label{lem:FundamentalLemmaDiracMultiplicationOperator}
	Let $\lambda_1,\lambda_2\ldots $ be the sequence of eigenvalues belonging to the smooth orthonormal basis of spinor fields $\varphi_1,\varphi_2,\ldots,$ of $L^2(S_F)$. 
	Let $k\in \set Z^2$ and $A(k)$ be the generalized multiplication operator in $l^2\oplus l^2=L^2(\set N, \set C^2)$ associated to the sequence $(M_l(k))_{l\in \set N}$ of complex $2\times 2$ matrices defined by
	\[
		M_{l}(k)=\left (\begin{matrix}- \lambda_l & \im(k_1-k_2) \\ \im(k_1+k_2) & \lambda_l \end{matrix} \right ).
	\]
	with 
	\[
		\dom{A(k)}=\{(x_l)_{l=1\ldots} \in l^2\oplus l^2| (M_l(k)x_l)_{l=1,\ldots} \in l^2\oplus l^2\}.
	\]
	
	Let $M_A$ be the operator of multiplication associated to $(A(k))_{k\in\set Z^n}$.
	Then for $\psi\in\Sec{S}$ holds:
	\[
		(\FT\circ L \circ B)(D\psi)=(M_A \circ \FT \circ L \circ B)\psi,
	\]
	or equivalently
	\[
		D\psi=(B^{-1}\circ L \circ \FT^{-1} \circ M_A \circ \FT \circ L \circ B)\psi.
	\]
\end{lemma}
\begin{proof}
			The $x^{\epsilon,k}_l:=\chi_{\{k\}}d^\epsilon_l$ for  $k\in \set Z^2, l \in \set N, \epsilon=\pm 1$ 
			form an orthonormal basis in $L^2(\set Z^2,l^2\oplus l^2)$, therefore it is sufficient to show
		\[
			\begin{split}
			\spr{(\FT\circ L \circ B)(D\psi),x^{\epsilon,k}_l}_{L^2(\set Z^2,l^2\oplus l^2)}=%
			&\spr{(M_A \circ \FT \circ L \circ B)(\psi),x^{\epsilon,k}_l}_{L^2(\set Z^2,l^2\oplus l^2)}\\
		&= \spr{(M_A\circ \FT \circ L \circ B)(\psi)(k),d^\epsilon_l}_{l^2\oplus l^2}\\
		&= \spr{ A(k) (\FT\circ L \circ B)(\psi)(k),d^\epsilon_l}_{l^2\oplus l^2}\\
	\end{split}
	\]
	for all $k,l,\epsilon$. 
	To this end, let $k,l,\epsilon$ be fixed. Then
		\[
			\begin{split}
				\spr{(\FT\circ L \circ B)(D\psi),x^{\epsilon,k}_l}_{L^2(\set Z^2,l^2\oplus l^2)}
&=\spr{(\FT\circ L \circ B)(D\psi)(k),d^\epsilon_l}_{l^2\oplus l^2}\\
&=\f{1}{(2\pi)^2} \int_{\set T^{1,1}} e^{-\im\spr{k,x}_2} \spr{( L \circ B)(D\psi),d^\epsilon_l}_{l^2\oplus l^2}\, dx\\
&=\f{1}{(2\pi)^2} \int_{\set T^{1,1}} e^{-\im\spr{k,x}_2} \spr{(B_\epsilon(D\psi))(x_1,x_2,\cdot),\varphi_l}_{L^2(S_F)} dx\\
&=\f{1}{(2\pi)^2} \int_{\set T^{1,1}} \int_F \sprBF{B_\epsilon(D(\psi))(x_1,x_2,\cdot), e^{\im\spr{k,x}_2}  \varphi_l} dx\\
\end{split}
\]
By lemma \ref{lem:SpinorBundleIsomorphism} (b) this equals to
\[
	=\f{1}{(2\pi)^2} \int_{\set T^{1,1}} \int_F \sprB{D \psi, \varphi^{\epsilon,k}_l}.
\]
By Fubini's theorem and because of $\sprB{\psi_1,\psi_2}=\isprB{\psi_1,J\psi_2}$ this in turn equals to
\[
	= \f{1}{(2\pi)^2} \int_{\set T^{1,1}\times F} \isprB{D \psi, J\varphi^{\epsilon,k}_l}
\]
and since $\im D$ is $J$-symmetric (\mcitesatz{3.18}{Baum1981}) this equals to 
\begin{equation}
	\label{eqn:pause}
	= \f{1}{(2\pi)^2} \int_{\set T^{1,1}\times F}- \isprB{ \psi,  D J \varphi^{\epsilon,k}_l}.
\end{equation}
Let us make a break here and just compute the integrand: $J\varphi^{\epsilon,k}_l=-\varphi^{-\epsilon,k}_l$ and by lemma \ref{lem:SpinorBundleIsomorphism} \ref{enum:IndefiniteSpinorBundle} 
	\[
		- \isprB{ \psi,  D J \varphi^{\epsilon,k}_l}=\isprB{\psi,D\varphi^{-\epsilon,k}_l}=-\sprBF{B_{-1}(\psi),B_1(D\varphi^{-\epsilon,k}_l)} - \sprBF{B_1(\psi),B_{-1}(D\varphi^{-\epsilon,k}_l)}.
 \]
 By lemma \ref{lem:DiracInSplittedBundle} this equals to 
 \[
	 -\sprBF{B_\epsilon(\psi),\epsilon e^{\im\spr{k,x}_2}\pr_2^*(D_F\varphi_l)} - \sprBF{B_{-\epsilon}(\psi),\im(k_1-\epsilon k_2) e^{\im\spr{k,x}_2}\pr_2^*(\varphi_l)}.
 \]
 Since $D_F\varphi_l=\lambda_l\varphi_l$ and $\lambda_l$ is real, this equals to 
	 \[
		 = - e^{-\im\spr{k,x}_2}\left ( \epsilon \lambda_l \sprBF{B_\epsilon(\psi),\pr_2^*(\varphi_l)}  +  \im(k_1-\epsilon k_2)\sprBF{B_{-\epsilon}(\psi),\pr_2^*(\varphi_l)} \right )
 \]
 Plugging this expression for the integrand into equation (\ref{eqn:pause}), we see that it equals to
	\[
		\begin{split}
			&=\f{1}{(2\pi)^2}\int_{\set T^{1,1}\times F} -e^{-\im\spr{k,x}_2}\left ( \epsilon \lambda_l \sprBF{B_\epsilon(\psi),\pr_2^*(\varphi_l)}  +  \im(k_1-\epsilon k_2)\sprBF{B_{-\epsilon}(\psi),\pr_2^*(\varphi_l)} \right )\\
			&=\f{1}{(2\pi)^2}\int_{\set T^{1,1}}-e^{-\im\spr{k,x}_2} \int_F \epsilon \lambda_l \sprBF{B_\epsilon (\psi),\varphi_l)}+\im(k_1-\epsilon k_2)\sprBF{B_{-\epsilon}(\psi),\varphi_l}\\
			&=\f{1}{(2\pi)^2}\int_{\set T^{1,1}}- e^{-\im\spr{k,x}_2} \left (\epsilon \lambda_l \spr{(L \circ B )(\psi),d^\epsilon_l}_{l^2\oplus l^2} +\im(k_1-\epsilon k_2)\spr{(L\circ B)(\psi),d^{-\epsilon}_l}_{l^2\oplus l^2} \right)\\
			&=- \epsilon \lambda_l \spr{(\FT \circ L \circ B)(\psi)(k),d^\epsilon_l}_{l^2\oplus l^2} + \im(k_1-\epsilon k_2)\spr{(\FT \circ L \circ B)(\psi)(k),d^{-\epsilon}_l}_{l^2\oplus l^2}\\
		&= \spr{ A(k) (\FT\circ L \circ B)(\psi)(k),d^\epsilon_l}_{l^2\oplus l^2}.
		\end{split}
	\]
\end{proof}

\begin{proposition}
	\label{prop:UnitaryEquivalenceOfDiracTorusProduct}
	Let $T:=\FT \circ L \circ B$. $T$ is a unitary equivalence of $L^2$ spaces. 
	$D$ is closable and for its closure $\bar D$ holds:
			\[
				\bar D \psi =(T^{-1} \circ M_A \circ T)\psi \quad \text{for }\psi \in \dom{\bar D}
\]
with
\[
	\dom{\bar D}= \{\psi \in L^2(S): T\psi \in \dom{M_A}\}.
\] 
\end{proposition}
\begin{proof}
	$T$ is a unitary equivalence since $\FT,L$ and $B$ are unitary equivalences.
	Let $M_0$ be the restriction of $M_A$ to 
	\[
		\dom{M_0}:=\{T \psi | \psi \in \Sec{S}\},
	\]
	such that $D=T^{-1}M_0T$ and $D$ and $M_0$ are unitarily equivalent.
	Since $A(k)$ is closed for all $k\in \set Z^2$ by proposition \ref{prop:GeneralizedMultiplicationOperatorOnZn} \ref{enum:MultiplicationOperatorClosed} and lemma \ref{lem:FundamentalLemmaDiracMultiplicationOperator}, the associated maximal generalized multiplication operator $M_A$ is closed by proposition \ref{prop:GeneralizedMultiplicationOperatorOnZn} \ref{enum:MultiplicationOperatorClosed}, and hence $M_0$ is closable and $\bar M_0\subset M_A$.
Furthermore, as $D$ and $M_0$ are unitarily equivalent (\ie $D=T^{-1}M_0T$ and $T$ is unitary), $D$ is closable and $\bar D=T^{-1}\bar M_0T$. We therefore need to prove that $\bar M_0=M_A$.

\begin{enumerate}
	\item Since $\bar M_0\subset M_A$ has already been shown, it remains to show that $M_A\subset \bar M_0$. Let $x^{\epsilon,k}_l$ be defined as in lemma 	\ref{lem:FundamentalLemmaDiracMultiplicationOperator}, then
\[
	\{x^{\epsilon,k}_l|\epsilon=\pm 1, k\in \set Z^2, l \in \set N\}\subset \dom{M_0},
\]
since $T\varphi^{\epsilon,l}_k=x^{\epsilon,l}_k$. Let $\abs k:=\max | k_1|,|k_2|$ for $k\in \set Z^2$.
Let now $x\in \dom {M_A}$ and let
\[
	x_j:=
	\sum_{1\leq l \leq j\atop {0 \leq |k|\leq j\atop\epsilon=\pm 1}} \spr{x(k),x^{\epsilon,k}_l}_{L^2(\set Z^2,l^2\oplus l^2)}x^{\epsilon,k}_l.
\]
$x_j\in L^2(\set Z^2,l^2\oplus l^2)$ is a finite linear combination of the $x^{\epsilon,k}_l$ and thus an element of $\dom{M_0}$. 
$x_j$ converges to $x$ \wrt the $L^2$-norm. $x_j$ also converges to $x$ pointwise, \ie for all $k\in \set Z^2$: $x_j(k)\to x(k)$.

\item \label{enum:PointwiseConvergence}
As a intermediate step we show for $k\in \set Z^2$ fixed 
\[
	A(k)x_j(k)\to A(k)x(k)
\]
in $l^2\oplus l^2$ for $j\to \infty$. By using the definition of $M_l(k)$, a simple calculation shows for $j>|k|$:
	\[
		\norm{A(k)x(k)-A(k)x_j(k)}^2_{l^2\oplus l^2}=\sum_{l>j,\atop\epsilon \pm 1} M_l(k) \left (\spr{x(k),d^1_l}_{l^2\oplus l^2} d_l + \spr{x(k),d^{-1}_l}_{l^2\oplus l^2} d^{-1}_l\right )
\]
and the last expression converges to $0$ by lemma since $x(k)\in \dom{A(k)}$.
\item 
Next we show $M_0 x_j \to M_A x$: Let now $\epsilon>0$, and choose $j_0\in \set N$ such that
\[
	\sum_{|k|>j_0} \norm{A(k)x(k)}^2_{l^2\oplus l^2} < \f{\epsilon}{2}.
\]
This is always possible since $M_Ax \in L^2(\set Z^2,l^2\oplus l^2)$. Let $\Omega_j:=\{k\in \set Z| |k|\leq j\}$. Now choose for all $k \in \set Z^2$ with $|k| \leq  j_0$  a $j_k\in \set N$ such that 
\[
	\norm{A(k)x(k)-A(k)x_{j'}(k)}^2_{l^2\oplus l^2}< \f{\epsilon}{2\# \Omega_{j_0}} \quad \text{ for all $j'>j_k$}.
\]
This is always possible since  $A(k)x_j(k)\to A(k)x(k)$ for all $k$ by \ref{enum:PointwiseConvergence}. Now let 
\[
	j:=\max \{j_k: |k|\leq j_0\} \cup \{j_0 \}
\]
Then for all $j'>j$ we have
\[
	\begin{split}
		\norm{M_A x - M_0 x_{j'}}&= \sum_{|k|>j'} \norm{A(k)x(k)}^2_{l^2\oplus l^2} + \sum_{|k|\leq {j'}}  	\norm{A(k)x(k)-A(k) x_{j'}(k)}^2_{l^2\oplus l^2}\\
		&< \sum_{|k|>{j_0}} \norm{A(k)x(k)}^2_{l^2\oplus l^2} +  \sum_{|k|\leq {j_0}}  	\norm{A(k)x(k)-A(k)  x_{j'}(k)}^2_{l^2\oplus l^2}\\
		&<\f{\epsilon}{2}+ \# \Omega_{j_0} \cdot \f{\epsilon}{2 \# \Omega_{j_0}}= \epsilon,
	\end{split}
\]
thus for all $x\in \dom {M_A}$ we have constructed a sequence  $(x_j)\in \dom {M_0}$ with $x_j\to x$ and $M_0 x_j \to M_A x$. This shows $x\in \dom{\bar M_0}$ and $\bar M_0 x=M_Ax$, hence $M_A \subset \bar M_0$.
\end{enumerate}
\end{proof}

We are now finally able to prove the main theorem of this section:

\begin{theorem}
	\label{thm:SpectrumTorusProduct}
	The spectrum of the Dirac operator $D$ of $\set T^{1,1}\times F$ is the whole complex plane $\set C$, in formulas:
	\[
		\sigma(D)=\sigma(\bar D)=\set C.
	\]
	The spectrum of $\bar D$ consists of two parts $\sigma(\bar D)=\specp(\bar D)\cup \specc(\bar D)$, the point spectrum $\specp(\bar D)$ and the continuous spectrum $\specc(\bar D)=\set C\setminus \specp(\bar D)$. 
	The residual spectrum $\specr(\bar D)$ is empty, and for the point spectrum holds
	\[
		\specp(D)=\specp(\bar D)=\{\pm\sqrt{\lambda_l^2-k_1^2+k_2^2}|l\in \set N, k \in \set Z^2\}.
	\]
\end{theorem}

\begin{proof}
	First note that $\sigma(D)=\sigma(\bar D)$ (\mcitep{105}{Mlak1991}). 
	Since $\bar D$ and $M_A$ are unitarily equivalent by proposition \ref{prop:UnitaryEquivalenceOfDiracTorusProduct}, it remains to calculate the spectrum of $M_A$.
	
	$A(k)$ is densely defined for all $k$ by proposition \ref{prop:GeneralizedMultiplicationOperatorOnZn} \ref{enum:ADenselyDefined}. 
	$\specr(A(k))=\emptyset $ for all $k$ by proposition \ref{prop:GeneralizedMultiplicationOperatorOnZn} \ref{enum:ResidualSpectrum}. Therefore we are in the scope of corollary \ref{cor:GeneralizedMultiplicationOperatorsOnZn} and immediately derive that the spectrum of $M_A$ consists only of the point spectrum $\specp(M_A)$ and the continuous spectrum $\specc(M_A)$:
	\[
		\sigma(M_A)=\specp(M_A)\cup\specc(M_A)	
	\]
with the following relations
\[
\specp(M_A)=\bigcup_{k \in \set Z^2} \specp(A(k))
\]
and
\[
	\specc(M_A)=\{\lambda \in \set C: \norm{(\lambda-A(k))^{-1}}\text{ is not bounded on } \set Z^2\}.
\]
%
Let us first determine the point spectrum of $M_A$: By applying proposition \ref{prop:GeneralizedMultiplicationOperatorOnZn} \ref{enum:PointSpectrum} twice, we find
\[
	\specp(M_A)=\bigcup_{k \in \set Z^2} \specp(A(k))=\bigcup_{k\in \set Z^2} \bigcup_{l\in \set N} \specp(M_l(k))
\]
	The spectrum of the $2\times 2$ matrix $M_l(k)$ consists of its eigenvalues, which are the zeroes of its characteristic polynomial 
	\[
		\lambda^2-\lambda_l^2+k_1^2-k_2^2,
	\]
	we therefore conclude
	\[
		\specp(M_A)=\{\pm \sqrt{\lambda_l^2-k_1^2+k_2^2}:  l \in \set N, k \in \set Z^2\}.
	\]	
	If $\lambda \in \set C \setminus \specp(M_A)$, then $\lambda-M_A$ is invertible and $\lambda-A(k)$ is invertible for all $k\in \set Z^2$ by (a). By lemma \ref{lem:InverseOfMultiplicationOperator} the inverse  of $\lambda-M_A$ is the multiplication operator associated to $(\lambda-A(k))^{-1}$.
	By proposition \ref{prop:GeneralizedMultiplicationOperatorOnZn} \ref{enum:MAInvertibleIFFAInvertible}, the inverse of $\lambda-A(k)$ is the generalized multiplication operator in $l^2\oplus l^2$ associated to the sequence 
	\[
		M_l'(k)=\f{1}{\lambda^2-\lambda_l^2+k_1^2-k_2^2}\left ( \begin{matrix} \lambda-\lambda_l & \im(k_1-k_2) \\ \im (k_1+k_2) & \lambda+\lambda_l \end{matrix} \right )
	\]
	Let us determine the continuous spectrum now: Let $a_1:=e_1-e_{2}, a_j:=j a_1$. 
	Then \[
		M_l'(a_j)=\f{1}{\lambda^2-\lambda_l^2} \left ( \begin{matrix} \lambda-\lambda_l & 2j \im \\ 0 & \lambda+\lambda_l\end{matrix} \right )
	\]
	Now the sequence $x_j$ defined by $x_j:=\chi_{\{a_j\}} d^{-1}_1$ fulfills $x_j\in \dom{(M_A)^{-1}}$, $\norm{x_j}=1$ and $\norm{(M_{\lambda-A})^{-1}x_j}^2=(\lambda+\lambda_l)^2 + 4j^2$, and since this last expression is unbounded this shows that  $(\lambda-M_A)^{-1}$ cannot be bounded on $\set Z^2$. 
	Therefore we conclude $\lambda \in \specc(M_A)$ by  corollary \ref{cor:GeneralizedMultiplicationOperatorsOnZn}.
\end{proof}

\begin{remark}
	\label{rem:TorusDifferentSpinStructures}
	In this remark we use theorem \ref{thm:SpectrumTorusProduct} to generalize the results of theorem \ref{thm:SpectrumTorus} to some (but not all) spin structures other than the $(0,\ldots,0)$ spin structure.
	We do not go into details, in particular we do not discuss the different spin structures.
	
	Note that we defined the torus $\set T^\pq$ to be isometric to $[0,2\pi]^{p+q}$. Other authors (in particular, \mcite{Ginoux2009} and \mcite{Lange}) defined the torus to be isometric to $[0,1]^{p+q}$ or used lattices to define the torus. In these cases, a normalization factor of $2\pi$ can occur. We adapted their results to our metric, such that the results are comparable with our theorems \ref{thm:SpectrumTorus} and \ref{thm:SpectrumTorusProduct}.


	The Dirac spectrum of the Riemannian $2N$-dimensional torus $\set T^{2N}$ was originally calculated by T. Friedrich (see \mcite{Ginoux2009} or the original article \mcite{FriedrichTorus}): Let $\delta=(\delta_1,\ldots,\delta_{2N})\in \{0,1\}^{2N}$.
	Then the spectrum of the Dirac operator $D_{\set T^{2N}_\delta}$ of $\set T^{2N}_\delta$ endowed with the flat metric and the $\delta$-spin structure is given by
	\[
		\specp(D_{\set T^{2N}_\delta})=\{\pm \sqrt{\sum_{j=1}^{2N} \left ( z_j + \f{1}{2} \delta_j z_j\right )^2}: z \in \set Z^{2N} \}
	\]
	Since the torus is compact and Riemannian, we have
	\[
		\specp(D_{\set T^{2N}_\delta})=\specp( \bar{D}_{\set T^{2N}_\delta}) =\sigma( \bar D_{\set T^{2N}_\delta}).
	\]
	Now $\set T^{1,1}\times \set T^{2N}_\delta=\set T^{1,2N+1}_{\delta'}$ carries the $\delta'$ spin structure with $\delta'=(0,0,\delta_1,\ldots,\delta_{2N})$ (thus $\delta'_j=0$ for $j=1,2$ and $\delta'_j=\delta_{j-2}$ for $j>2$).
	Let $D$ be its Dirac operator.
	Using theorem \ref{thm:SpectrumTorusProduct} we see that 
	\[
		\sigma(D)=\sigma(\bar D)=\specp(\bar D)\cup \specc(\bar D)
	\]
	with
	\[
		\begin{split}
		\specp(\bar D)&=\{\pm  \sqrt{ \sqrt{\sum_{j=1}^{2N} \left ( z_j + \f{1}{2} \delta_j z_j\right )^2 } -k_1^2+ k_2^2} : z\in \set Z^{2N}, k \in \set Z^2\}\\
		&=\{\pm \sqrt{\sum_{j=1}^{2N+2} \kappa_j\left ( (k_j+\f{1}{2} \delta'_j k_j \right )^2}: k \in \set Z^{2N+2}\}
	\end{split}
	\]	
and 
\[
	\specc(\bar D)=\set C\setminus \specp(\bar D).
\]
	Note that the point spectrum of $\set T^\pq_\delta$ with arbitrary spin structure $\delta=(\delta_1,\ldots,\delta_n)\in \{0,1\}^n$ has been determined before, see \mcitesatz{4.1.2}{Lange}).
\end{remark}

\bibliographystyle{alpha}
\bibliography{../bibliography-diplom.bib}{}
\end{document}